\documentclass[11pt]{amsart}
    \usepackage{amsmath,amsfonts,eucal,xy, enumerate, amscd}

\xyoption{all}

     \newtheorem{theorem}{Theorem}[section]
     \newtheorem{proposition}[theorem]{Proposition}

     \newtheorem{lemma}[theorem]{Lemma}

 \theoremstyle{definition}

     \newtheorem{example}[theorem]{Example}
     
     \newtheorem{definition}[theorem]{Definition}
     \newtheorem{remark}[theorem]{Remark}

     \newcommand{\In}{\subseteq} 
\newcommand{\script}[1]{\mathcal #1}
\renewcommand{\a}{\alpha}
\newcommand{\as}{{\operatorname{as}}}
\renewcommand{\b}{\beta}
\renewcommand{\Bbb}[1]{\mathbb #1}
\newcommand{\bb}{{\Bbb B}}
\newcommand{\bc}{{\Bbb C}}
\newcommand{\bk}{{\Bbb K}}
\newcommand{\bl}{{\Bbb L}}

\newcommand{\br}{{\Bbb R}}
\newcommand{\bsur}{{\operatorname{bsur}}}

\newcommand{\conv}{{\operatorname{conv}}}

\newcommand{\cM}{{\mathcal M}}
\newcommand{\cP}{{\mathcal P}}

\newcommand{\dual}{{\operatorname{dual}}}

\newcommand{\fin}{{\operatorname{fin}}}
\newcommand{\g}{\gamma}

\newcommand{\Groth}{{\operatorname{Groth}}}

\renewcommand{\frak}{\mathfrak}
\newcommand{\gA}{{\frak A}}

\newcommand{\gF}{{\frak F}}
\newcommand{\gK}{{\frak K}}
\newcommand{\gL}{{\frak L}}
\newcommand{\gM}{{\frak M}}
\newcommand{\gN}{{\frak N}}
\newcommand{\gP}{{\frak P}}
\newcommand{\gW}{{\frak W}}
\newcommand{\inj}{{\operatorname{inj}}}
\newcommand{\lam}{\lambda}

\newcommand{\loc}{{\operatorname{loc}}}
\newcommand{\lup}{{\operatorname{lup}}}

\newcommand{\ori}{{\operatorname{ori}}}

\newcommand{\po}{^\circ}

\newcommand{\resp}{\rm {resp.\,}}
\newcommand{\rup}{{\operatorname{rup}}}
\newcommand{\s}{\sigma}
\newcommand{\sB}{{\script B}}

\newcommand{\sC}{{\script C}}

\newcommand{\sG}{\cP}

\newcommand{\sL}{{\script L}}
\newcommand{\sM}{{\script M}}
\newcommand{\sO}{{\script O}}
\newcommand{\sT}{{\script T}}

\newcommand{\sur}{{\operatorname{sur}}}

\newcommand{\V}{\Vert}
\newcommand{\von}{{\operatorname{von}}}

\newcommand{\wF}{\widetilde F}

\newcommand{\wQ}{\widetilde Q}

\newcommand{\wX}{\widetilde X}
\newcommand{\wY}{\widetilde Y}

     \newcommand{\ie}{{i.e., }}
     \newcommand{\eg}{{e.g.\  }}
    \newcommand{\cf}{{cf.\ }}

     \newcommand{\Sset}{\supseteq}
       
     \newcommand{\lrw}{\longrightarrow}

\begin{document}

\parskip=0.5\baselineskip
\baselineskip=1.44\baselineskip

\title
{The triangle of operators, topologies, bornologies}
\date{May 31, 2005}
\author{Ngai-Ching Wong}

\address{
Department of Applied Mathematics, National Sun Yat-sen
  University, and National Center for Theoretical Sciences, Kaohsiung, 80424, Taiwan, R.O.C.}

\email{wong@math.nsysu.edu.tw}

\dedicatory{In memory of my teacher Yau-Chuen Wong
(1935.10.2--1994.11.7)}

\thanks{Partially supported  by
Taiwan National Science Council, NSC
83-0208--M--110--017, 93-2115-M-110-009}



\begin{abstract}
This paper discusses two common techniques in functional analysis:
the topological method and the bornological method. In terms of
Pietsch's operator ideals, we establish the equivalence of the
notions of operators, topologies and bornologies. The
approaches in the study of locally convex spaces of Grothendieck
(via Banach space operators), Randtke (via continuous seminorms) and
Hogbe-Nlend (via convex bounded sets) are compared.
\end{abstract}

\maketitle

\tableofcontents

\section{
Introduction }

How can one describe  a linear operator $T$ from a Banach space $E$
into a Banach space $F$? The usual way to describe $T$ is to state
either  the bornological property, via $TU_E$, or  the topological
property, via $T^{-1}U_F$, of $T$, where $U_E$ (resp.\ $U_F$) is the
closed unit ball of $E$ (resp.\ $F$).  However, there are a lot of
examples indicating that these two machineries are equivalent. For
instance,
\begin{itemize}
    \item \quad\ $T$ is {\em bounded} (\ie $TU_E$ is a bounded subset of
$F$) \\ $\Leftrightarrow$ $T$ is {\em continuous} (\ie $T^{-1}U_F$ is a
0-neighborhood of $E$ in the norm topology);
    \item \quad\ $T$ is {\em of finite
rank} (\ie $TU_E \subseteq \conv\{y_1, y_2, \ldots, y_n\}$ for some
$y_1$, $y_2$, $\ldots$, $y_n$ in $F$) \\ $\Leftrightarrow$ $T$ is {\em
weak-norm continuous} (\ie $T^{-1}U_F$ is a 0-neighborhood of $E$ in
the weak topology); and
    \item \quad\ $T$ is {\em compact} (\ie $TU_E$ is totally
bounded in
 $F$) \\  $\Leftrightarrow$ $T$ is continuous in the topology of uniform
 convergence on norm compact subsets of $E'$
 (\ie $T^{-1}U_F \supseteq K^\circ$, the polar of a norm compact  subset $K$ of
   the dual space $E'$ of $E$).
\end{itemize}
This is because the unit ball of a normed space simultaneously
serves as a neighborhood of zero and a bounded set.   It is,
however, no longer true in the context of locally convex spaces
(LCS's, shortly).  Mackey-Arens' Theorem indicates that topologies
(families of neighborhoods) and bornologies (families of bounded
sets) are in dual pair (see \eg \cite{Sch71}).

It is a long tradition of classifying special classes of
locally convex spaces  by families of continuous
operators among them. A famous example is, of course,
 Grothendieck's identification of the  class of nuclear locally convex spaces.
Other examples are those of Schwartz LCS's, infra--Schwartz
LCS's and their ``co--spaces". After the great effort of Pietsch
\cite{Pie80}, it is now well--known that such suitable families of
continuous operators are the so--called operator ideals.

There are many ways to utilize Grothendieck's idea. For example, one
can
define a LCS $X$ to be nuclear ($\resp$ Schwartz, infra--Schwartz) by
asking
that for each continuous seminorm $p$ on $X$, there is a continuous
seminorm $q$ on $X$ with $p\leq q$ such that the canonical map
$\wQ_{pq}$ from $\wX_q
=\widetilde{X/q^{-1}(0)}$ into $\wX_p=\widetilde{X/p^{-1}(0)}$ is
nuclear
($\resp$ precompact,
weakly compact), where \ $\widetilde{\vphantom{I}}$ \ denotes
completion. It
amounts to saying that the completion $\wX$ of $X$ is a topological
projective
limit $\varprojlim\wQ_{pq}\wX_q$ of Banach spaces of nuclear
type ($\resp$ precompact type, weakly compact type). The converse is
also true, see Junek \cite[p.~139]{Jun83}.
We call such a LCS a Grothendieck space of nuclear ($\resp$
precompact, weakly compact) type, or shortly a $\Groth (\gN)$--space
($\resp$ $\Groth(\gK_p)$--space, $\Groth(\gW)$--space), where $\gN\ (\resp\gK_p,\gW)$
is the
ideal of all nuclear ($\resp$ precompact, weakly compact) operators
between
\emph{Banach spaces}.

 As a
dual
concept, a locally convex space $X$ is said to be  a co--Grothendieck space of type $\gA$, or  shortly a
co--$\Groth(\gA)$--space,
if for each infracomplete
disk $A$
in $X$ there is an infracomplete disk $B$ in $X$ such that $A\subseteq B$
and
the canonical map $J_{BA}$ from
$X(A)=\bigcup_{\lam>0}\lam A$ into
$X(B)=\bigcup_{\lam>0}\lam B$
belongs to $\gA(X(A),X(B))$.
In other words, the convex bornological vector space $X$ equipped with
the
infracomplete bornology of $X$ is the bornological inductive limit
$\varinjlim
J_{BA}X(A)$ of Banach spaces of type $\gA$. The converse is again true.

Another way to go is to define the ideal topology and the ideal
bornology on
each LCS associated to an operator ideal $\gA$ \emph{on LCS's}. A
continuous
seminorm $p$ on a LCS $X$ is said to be an $\gA$--continuous seminorm
if the
canonical map $\wQ_p :X\to\wX_p$ belongs to the injective hull
$\gA^{\inj}$ of
$\gA$. The topology on $X$ defined by the family of all such seminorms
is
called the $\gA$--topology of $X$.
Similarly, an absolutely convex bounded set $B$ in $X$
is said to be $\gA$--bounded if the canonical map $J_{B}$ from $X(B)=\bigcup_{\lam>0}\lam B$
into $X$ belongs to the bornologically surjective hull $\gA^{\bsur}$ of
$\gA$.  The bornology on $X$
defined by the family of all such bounded sets is called the
$\gA$--bornology of $X$.
A LCS $X$ is said to be $\gA$--topological (resp. $\gA$--bornological)
if the topology (resp. bornology) of $X$ coincides with the
$\gA$--topology (resp. $\gA$--bornology).

In \cite{W94} we show that Grothendieck spaces are
essentially a kind of $\gA$--spaces.  Thus these two different approaches
coincide.  In this paper, we will develop the duality theory of
$\gA$--topological spaces
and $\gA$--bornological spaces. Basically, one may expect that a
locally convex
space $X$ is $\gA$--topological ($\resp\gA$--bornological) if and only
if its
strong dual $X'_\b$ is $\gA$--bornological ($\resp\gA$--topological).
One can discover the same is true for
Grothendieck spaces and co--Grothendieck spaces by observing
 the duality of topology and bornology and
the duality
of projective limits and inductive limits (see, e.g.,  \cite{Jar81, Jun83}).

The  following
commutative diagram summaries our works.

$$
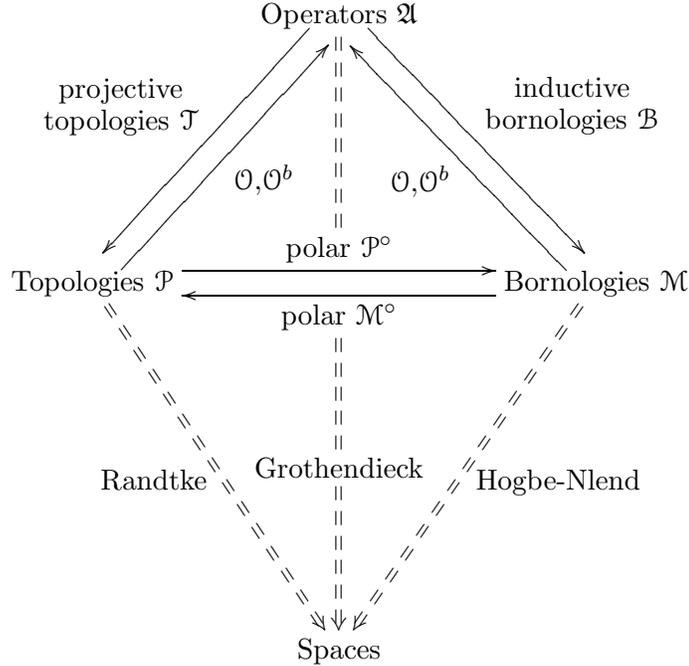
\xymatrix{ & \txt{Operators $\gA$}
    \ar @<-1ex> [dddl]_-{\txt{projective\\topologies $\sT$}}
    \ar @<1ex> [dddr]^-{\txt{inductive\\bornologies $\sB$}}
\\
\\
\\
\txt{Topologies $\cP$}
    \ar @<-1ex> [uuur]_-{\txt{$\sO$,$\sO^b$}}
    \ar @<1ex> [rr]^-{\txt{polar $\cP^\circ$}}="t"
    \ar @2{-->} [ddddr]_-{\txt{Randtke}}
&& \txt{Bornologies $\cM$}
    \ar @<1ex> [uuul]^-{\txt{$\sO$,$\sO^b$}}
    \ar @<1ex> [ll]^-{\txt{polar $\cM^\circ$}}="b"
    \ar @2{-->} [ddddl]^-{\txt{Hogbe-Nlend}}
\\
\\
& {\txt{Grothendieck}}
\\
\ar @2{--} "6,2";"b"
\\
& {\txt{Spaces}} \ar @2{--} "1,2";"t" \ar @2{-->} "6,2";"8,2" }
$$
The theory of operator ideal is founded by Pietsch \cite{Pie80} and
originated from the works of Grothendieck \cite{Gro56} and Schatten
\cite{Schatten60}.  See also \cite{Mic78, Jar81, Jun83, Def93} for
more information.  The idea of generating topologies and generating
bornologies are due to Stephani \cite{Step70, Step72, Step73,
Step80, Step83} and Franco and Pi\~neiro \cite{FraPin82} in the
context of Banach spaces. The explicit construction (with all arrows
shown in the diagram) of the (upper) triangle is given in \cite{WW88}, in
which several applications to Banach space theory are demonstrated.
When the underlying space is a {\em fixed} complex Hilbert space,
West implements the triangle in the context of operator algebras
\cite{West} and provides several applications with Conradie
\cite{CanWest} (see Section \ref{s:eg-HB}).
In
this paper, we shall complete the LCS version of the triangle.
As an application, we  shall show that in the study of LCS's,
the topological machinery of Randtke (via continuous
seminorms) \cite{Ran72} or the bornological machinery of Hogbe-Nlend
(via convex bounded subsets) \cite{HN77, HN81} is as strong as that
of the operator theoretical machinery of Grothendieck (via Banach
space operators) (see \eg \cite{Wong76, Wong79, Jun83, W94}).

The author dedicates this paper to his late teacher, Professor Yau-Chuen Wong, who
 introduced the same concept of $\gA$--topology and $\gA$-bornology through
a
great number of examples of special LCS's as well as partially
ordered locally convex
spaces (see, \cite{{Wong76},{Wong79}}), although he did not employ the Pietsch's language (operator
ideals) at
his time.
Together with \cite{WW88, W94}, the current paper is a continuation of his ideal (see \cite{W96}).

\section{Established examples in Hilbert spaces and Banach spaces}\label{s:eg-HB}

\subsection{The triangle for Hilbert spaces}\label{ss:tri_H}

Let ${\gM}$ be a von Neumann algebra of bounded linear operators on
a Hilbert space $H$, and ${\gA}$ an arbitrary non-zero
two-sided  ideal of ${\gM}$.

\begin{itemize}
    \item A locally
convex topology $\cP$ of $H$ is called a {\em generating topology}
if $\cP$ consists of norm open sets in $H$ such that all
operators in ${\gM}$ are $\cP$-to-$\cP$ continuous on $H$, \ie
the pre-images of $\cP$-open sets being $\cP$-open.

    \item A convex
vector bornology $\cM$ of $H$ is called a {\em generating bornology}
if $\cM$ consists of norm bounded subsets of $H$ such that all
operators in ${\gM}$ are $\cM$-to-$\cM$ bounded, \ie  sending
$\cM$-bounded sets to $\cM$-bounded sets.

    \item The $\gA$-\emph{topology} $\sT({\gA})$ is the projective topology of $H$ induced by operators in ${\gA}$, \ie
the weakest locally convex topology $t$ of $H$ such that
operators in ${\gA}$ are $t$-to-norm continuous.

    \item The $\gA$-\emph{bornology} $\sB({\gA})$ is the inductive bornology of $H$ induced by
operators in ${\gA}$, \ie the smallest convex vector bornology $b$
of $H$ such that operators in ${\gA}$ are norm-to-$b$ bounded.

    \item The polar of a subset $A$ in $H$ is
    $$
    A\po=\{x\in H: |\left<a,x\right>|\leq 1, \forall a\in A\}.
    $$

\end{itemize}

Remark that the ideal
${\gA}$ is
\begin{itemize}
    \item {\emph{self-adjoint}}, \ie  $T \in {\gA}$ if and only if its
Hilbert space adjoint map $T^*\in {\gA}$;
    \item {\em injective}, \ie  $T \in {\gA}$ whenever $\|Th\| \leq \|Sh\|, \forall h \in H$,
    for any $S$ in ${\gA}$ and $T$ in ${\gM}$; and
    \item {\em surjective}, \ie $T \in {\gA}$ whenever $TU_H \subseteq SU_H$
for any $S$ in ${\gA}$ and $T$ in ${\gM}$.
\end{itemize}

\begin{theorem}[{West \cite{West}}]\label{thm:west}
\
\begin{enumerate} 
    \item \begin{enumerate}
        \item The $\gA$-topology $\sT({\gA})$  is a generating
topology.
        \item  The $\gA$-bornology $\sB({\gA})$  is a
generating bornology.
    \end{enumerate}

        \item \begin{enumerate}
            \item The set $\sO(\cP)=\gL(H_\cP,H)$ of all $\cP$-to-norm continuous
        linear operators on $H$ is a  two-sided ideals of ${\gM}$.
            \item The set $\sO(\cM)=B(H,H^\cM)$ of all norm-to-$\cM$ bounded
        linear operators on $H$ is a two-sided ideals of ${\gM}$.
        \end{enumerate}
        \item \begin{enumerate}
            \item The polar $\cP\po=\{B\In H: B\po\text{ is a }\cP\text{-neighborhood of }0\}$
        of a generating topology $\cP$ is a generating bornology.
            \item The polar $\cM\po=\{V\In H: V\po\text{ is }\cM\text{-bounded}\}$ of a generating bornology $\cM$ is a generating topology.
            \end{enumerate}

    \item The triangle of operators, topologies and bornologies is commutative:
        \begin{enumerate} 
            \item $\sO(\sT({\gA})) = {\gA}$, $\sO(\sB({\gA})) = {\gA}$.
            \item $\sT({\gA})\po = \sB({\gA})$, $\sB({\gA})\po = \sT({\gA})$.
            \item $\sT(\sO(\cP)) = \cP$, $\sB(\sO(\cM)) = \cM$.
        \end{enumerate}
\end{enumerate}
\end{theorem}

\subsection{The triangle for Banach spaces}\label{ss:tri_B}

The   Banach space version of the ``triangle''  is known to have
many applications (\cf \cite{Pie80, Jar81, Jun83}). Let
$\gA=\bigcup\{\gA(E,F): E, F \text{ are Banach spaces}\}$ be an
{\em operator ideal} on Banach spaces in the sense of Pietsch
\cite{Pie80}:
\begin{enumerate}
    \item[({$\text{OI}_1$})]  The components $\gA(E,F)$ of $\gA$ are non-zero
subspaces of $\gL(E,F)$.
    \item[({$\text{OI}_2$})]  $RTS \in \gA(E_0,F_0)$ whenever $R
\in \gL(F,F_0)$, $T \in \gA(E,F)$ and $S \in \gL(E_0,E)$ for
arbitrary Banach spaces $E, E_0, F$, and $F_0$.
\end{enumerate}
An operator ideal $\gA$ is said to be {\em symmetric} if $T \in \gA(E,F)$ ensures
its Banach space dual map $T' \in \gA(F',E')$, and
{\em completely symmetric} if $T \in
\gA(E,F) \Leftrightarrow T' \in \gA(F',E')$.

Suppose for each Banach space $E$, we have
a  locally convex topology $\cP(E)$  consisting of norm
open subsets of $E$ and  a  convex vector bornology $\cM(E)$  consisting of norm
bounded subsets of $E$.
We call $\cP = \{\cP(E): E \text{ is a Banach space}\}$ a generating topology and
$\cM = \{\cM(E) : E \text{ is a Banach space}\}$
a generating bornology on Banach spaces, if
operators in $\gL(E,F)$ are $\cP(E)$-to-$\cP(F)$
continuous and $\cM(E)$-to-$\cM(F)$ bounded for all Banach spaces
$E$ and $F$, respectively.

The polar $\cP\po$ of a generating topology $\cP$ consists of components
$$
\cP\po(E)=\{B\in E: B\po \text{ is } \cP(E')\text{-bounded}\}.
$$
Similarly, the polar $\cM\po$ of a generating bornology $\cM$ consists of components
$$
\cM\po(E)=\{V\in E: V\po \text{ is a } \cP(E')\text{-neighborhood of } 0\}.
$$

\begin{theorem}[{Stephani  \cite{Step70, Step72, Step73,
Step80, Step83} and Wong and Wong \cite{WW88}}]\label{thm:Bspaceversion}
\
\begin{enumerate} 
    \item \begin{enumerate}
        \item The family of projective topologies $\sT(\gA)(E)$ of  Banach spaces
    $E$ induced by operators in $\gA(E,\,\cdot\,)$  forms a generating topology.
        \item The family of inductive bornologies $\sB(\gA)(F)$ of  Banach spaces $F$
induced by operators in $\gA(\,\cdot\,,F)$ forms a generating
bornology.
    \end{enumerate}

    \item \begin{enumerate}
        \item The family of sets $\sO(\cP)(E,F)=\gL(E_\cP,F)$ of all $\cP(E)$-to-norm continuous
    linear operators from $E$ into $F$ forms an injective operator ideal.
        \item The family of sets $\sO(\cM)(E,F)
=B(E,F_\cM)$ of all norm-to-$\cM(F)$ bounded linear operators from
$E$ into $F$
 forms a surjective operator ideal.
    \end{enumerate}

    \item \begin{enumerate}
        \item The polar $\cP\po$ of a generating topology $\cP$
is a generating bornology.
        \item The polar $\cM\po$ of a generating
bornology $\cM$ is a generating topology.
    \end{enumerate}

    \item The triangle of operators, topologies and bornologies is almost commutative.
    \begin{enumerate} 
            \item $\sO(\sT(\gA)) = \gA^\inj$
            and $\sO(\sB(\gA)) = \gA^\sur$.
            \item $\sB(\gA)\po = \sT(\gA)$ if\/ $\gA$ is {symmetric}, and $\sT(\gA)\po =
\sB(\gA)$ if\/ $\gA$ is {completely symmetric}.
            \item $\sT(\sO(\cP)) = \cP$, $\sB(\sO(\cM)) = \cM$.
    \end{enumerate}
\end{enumerate}
\end{theorem}

Note that the ideals $\gL$ of all bounded operators, $\gF$ of all
bounded operator of finite rank and $\gK$ of all compact operators
are all injective, surjective and completely symmetric.  These
explain the equivalence of topological and bornological approaches
for these operators demonstrated at the very beginning of this paper.
On the other hand, the ideal
$\gN$ of nuclear operators is neither injective, surjective or
completely symmetric (\cf \cite{Pie80}).


\section{Notations and Preliminaries}\label{s:preliminary}

The classic reference to the theory of operator ideals is, of course,
Pietsch
\cite{Pie80}. See also Jarchow \cite{Jar81} and Junek \cite{Jun83}. For the theory of
locally
convex spaces, together with Wong \cite{Wong92},
Schaefer \cite{Sch71} is our favorite. Hogbe--Nlend \cite{HN77} serves as
our
main source of the theory of bornology.

Throughout this paper, all vector spaces have the same underlying
scalar field $\bk$.
$\bk$ is either the field $\br$ of real numbers or the field $\bc$
of
complex numbers.
Locally convex topologies are always Hausdorff, and
convex vector bornologies are always separated, \ie
no nonzero subspace is bounded. Operators always refer
to
linear maps without any topological or bornological assumption.
$U_N$ always denotes the {\it closed\/} unit ball of a normed space
$N$.

A subset $B$ of a LCS $X$ is said to be a {\it disk\/} if $B$ is
{\it absolutely convex,\/} \ie $\alpha B+\b B\subseteq B$ whenever
$|\alpha|+|\b|
\leq1$. A disk $B$ is said to be a $\s$--{\it disk,\/} or {\it
absolutely\/}
$\s$--{\it convex\/} if $\Sigma_{n}\lam_n b_n$ converges in
$B$ whenever $\sum_n |\lam_n| \leq 1$ and $b_n\in B$,
$n=1,2,
\dots$. A bounded disk $B$ is said to be {\it infracomplete\/} (or a
{\it Banach disk\/}) if the normed space
$X(B)=\bigcup_{\lam>0} \lam B$
equipped with the gauge $\g_B$ of $B$ as its norm is complete, where
$\g_B(x)=
\inf\{\lam > 0 :x\in\lam B\}$, for each $x$ in $X(B)$. Any continuous
image of a
$\s$--disk or an infracomplete bounded disk is still a $\s$--disk or an
infracomplete bounded disk, respectively.
A LCS $X$ is said to be {\it infracomplete\/} if the {\it von Neumann
bornology\/}
$\cM_{\von}(X)$, \ie the original bornology induced by the topology
of $X$,
has a basis consisting of infracomplete subsets of $X$, or
equivalently,
$\s$--disked subsets of $X$.
In other words, $(X,\cM_{\von}(X))$ is a complete convex bornological
vector space.

Let $\langle X,X'\rangle$ be a dual pair and $B\subseteq X$. The
(absolute) polar
$B^\circ$ of $B$ in $X'$ is defined by
$$
B^{\circ}=\{x\in X' :|\langle b,x\rangle|\leq1,\quad\forall b\in B\}.
$$
Whenever $A\subseteq X'$,   denote by $A^\bullet$ the polar of $A$ taken
in
$X''_{\b\b}$, namely,
$$
A^\bullet=\{x\in X''_{\b\b} :|\langle a,x\rangle|\leq1,\quad\forall
a\in A\},
$$
where $X''_{\b\b}$ is the {\it strong bidual\/} of $X$, while $A^\circ$
denotes the polar of $A$ taken in $X$ with respect to the dual pair $\langle X,X'\rangle$.

\begin{proposition}[{See, e.g., Wong \cite[{pp.\ 224 and 227}]{Wong92}}]\label{prop:2.1}
Let $X$ and $Y$ be LCS's and
$T\in\sL(X,Y)$. We
have
\begin{enumerate}[(1)]
    \item $T\in\sL(X_\s,Y_\s)$, where $X_\s$, $Y_\s$ denote the
LCS's
in their weak topologies.
    \item $T\in\sL(X_\tau,Y_\tau)$, where $X_\tau$, $Y_\tau$ denote
the
LCS's in their Mackey topologies.
    \item $T'\in\sL(Y'_\b,X'_\b)$, where $T'$ is the dual map of $T$
and
$X'_\b$ (resp.\ $Y'_\b$) is the strong dual of $X$ (resp.\ $Y$).
    \item $(TA)^{\circ}=(T')^{-1}A^{\circ}$ for all nonempty subset
$A$ of $X$.
    \item $(T'B)^{\circ}=T^{-1}B^{\circ}$ for all nonempty subset
$B$ of $Y'$.
    \item $(T^{-1}W)^{\circ}=T'W^{\circ}$ for all
neighborhoods $W$ of $0$ in its Mackey topology $\tau(Y,Y')$.
\end{enumerate}
\end{proposition}

Let $X$ and $Y$ be LCS's.
$J$ in $\gL(X,Y)$ is called a {\em (topological) injection} if $J$ is one-to-one and relatively open.
$Q$ in $\gL(X,Y)$ is called a {\em (topological) surjection} if $Q$ is open (and thus $Q$ induces the topology of
$Y$).
$Q^1$ in $\gL(X,Y)$ is called a
{\it bornological surjection\/} if $Q^1$ is onto and induces the
bornology of
$Y$ (\ie for each bounded subset $B$ of $Y$ there is a bounded subset
$A$ of
$X$ such that $Q^1A=B$).

An operator ideal $\gA$ on LCS's is said to be
\begin{itemize}
    \item \emph{injective\/} if $JT\in\gA(X,Y_0)$ infers
$T\in \gA(X,Y)$, whenever  $T\in \gL(X,Y)$ and
$J\in\gL(Y,Y_0)$ is an injection for some LCS $Y_0$;
    \item \emph{surjective\/} if $TQ\in\gA(X_0, Y)$ infers
$T\in \gA(X,Y)$, whenever $T\in \gL(X,Y)$ and
$Q\in\gL(X_0, X)$  is a surjection  for some LCS $X_0$; and
    \item \emph{bornologically
surjective\/} if $TQ^1 \in \gA(X_0,Y)$ infers
$T\in \gA(X,Y)$,  whenever  $T\in \gL(X,Y)$ and
$Q^1\in\gL(X_0,X)$  is a bornological surjection
for some LCS $X_0$.
\end{itemize}

The \emph{injective hull} $\gA^\inj$, the \emph{surjective hull} $\gA^\sur$,
 and the \emph{bornologically surjective hull} $\gA^{\bsur}$ of $\gA$ is the
intersection of all injective, surjective, and bornologically surjective operator ideals
containing $\gA$, respectively.
Note that for operator ideals on Banach spaces, the notions of surjectivity and bornological surjectivity coincide.

Associate to each normed space $N$ the Banach space
$N^{\inj}=\l_\infty(U_{N'})$ and
the injection $J_N$ in $\gL(N,N^{\inj})$ defined by
$J_N(x)={(<x,a>)}_{a \in U_{N'}}$.
Similarly, we define $N^{\sur}$ to be the
{\it normed\/} space $L_1(U_N)=\{(\lambda_x)\in\ell_1(U_N) :
\sum_{x\in U_N}\lam_x\,x$ converges in  $N\}$
and $Q_N :N^{\sur}\to N$ to be the
surjection defined by $Q_N((\lam_x)_{x\in U_N})=\sum_{x\in U_N}
\lam_x\, x$.
In case $E$ is a Banach space, it is well--known that $E^{\inj}$ has
the  extension property and
$E^{\sur}$ has the  lifting property, cf. \cite{Pie80}.

\begin{proposition}[\cite{Pie80, FraPin82, WW93}]\label{hulls}
\
\begin{enumerate} 
    \item Let $\gA$ be an operator ideal on Banach spaces.
    \begin{eqnarray*}
        \gA^{\inj}(E,F)&=&\{R\in\gL(E,F) : J_FR\in\gA(E,F^{\inj})\}, \\
         \gA^{\sur}(E,F)&=&\{S\in\gL(E,F) : SQ_E\in\gA(E^{\sur},F)\}.
    \end{eqnarray*}

    \item Let $\gA$ be an operator ideal on LCS's.
        We can associate to each LCS $Y$ a LCS $Y^\infty$ and an
injection $J^\infty_Y$ from $Y$ into $Y^\infty$, and to each LCS  $X$ a LCS $X^1$ and a
bornological surjection $Q^1_X$ from $X^1$ onto $X$ such that
    \begin{eqnarray*}
\gA^{\inj}(X,Y)&=&\{R\in\gL(X,Y) :J^\infty_YR\in\gA(X,Y^\infty)\},\\
\gA^{\bsur}(X,Y)&=&\{S\in\gL(X,Y) :SQ^1_X\in\gA(X,Y^1)\}.
    \end{eqnarray*}
Moreover, we have
$$
\gA^{\inj\,\bsur} = \gA^{\bsur\,\inj}.
$$
In case $N$ is a normed space, $R\in\gL(X,N)$ and $S\in\gL(N,Y)$,
    \begin{eqnarray*}
J_NR\in\gA(X,N^{\inj}) &\iff& J^\infty_NR\in\gA(X,N^\infty), \\
SQ^1_N\in\gA(N^1,Y) &\iff& SQ_N\in\gA(N^{\sur},Y).
    \end{eqnarray*}

    \end{enumerate}
\end{proposition}


\section{The construction and the commutativity of the triangle}\label{s:triangle}

 Let $\sC$ be a  class of locally convex spaces.
 Let $X, Y\in\sC$. We denote by $\gL^b(X,Y)$, $\gL(X,Y)$ and
$L^{\times}(X,Y)$ the collection of all operators from $X$ into $Y$
which are bounded (\ie sending a 0-neighborhood to a bounded set),
continuous, and locally bounded (\ie sending bounded sets to bounded
sets), respectively.

Denote by $\sigma(X,X')$  the weak topology
of $X$ with respect to its dual space $X'$, while  $\cP_\ori(X)$ is the original topology of $X$.
We employ the notion $\cM_\fin(Y)$ for the {\em finite dimensional bornology} of $Y$ which has a basis consisting of all convex hulls
of finite sets.
On the other hand, $\cM_\von(Y)$ is used for the {\em von Neumann bornology} of $Y$
which consists of all topologically bounded subsets of $Y$.
Ordering of topologies and bornologies are induced by set-theoretical inclusion, as usual.  Moreover, we write briefly
$X_{\cP}$ for a vector space $X$ equipped with a locally convex
topology $\cP$ and
$Y^{\cM}$ for a vector space $Y$ equipped with a convex vector bornology
$\cM$.

We now give the details of the ``triangle''.

\begin{definition}\
\begin{enumerate}
    \item ({\bf ``Operators''}) A family $\gA = \{\gA(X,Y) : X,Y \in \sC\}$ of algebras of operators associated to each pair of spaces $X$ and $Y$ in $\sC$ is called an {\em operator ideal}
if
    \begin{description}
        \item[OI$_1$]  $\gA(X,Y)$ is a nonzero vector subspace of $\gL(X,Y)$ for all $X$, $Y$ in $\sC$; and
        \item[OI$_2$]  $RTS \in \gA(X_0,Y_0)$ whenever $R \in \gL(Y,Y_0)$,
$T\in\gA(X,Y)$ and $S \in \gL(X_0,X)$ for any $X_0$, $X$, $Y$ and $Y_0$ in $\sC$.
    \end{description}

    \item ({\bf ``Topologies''}) A family $\cP =\{\cP(X) : X\in \sC\}$ of locally convex   topologies
associated to each space $X$ in $\sC$ is called a {\em generating topology} if
    \begin{description}
    \item[GT$_1$] $\sigma(X,X') \In \cP(X) \In \cP_\ori(X)$ for all $X$ in $\sC$; and
    \item[GT$_2$] $\gL(X,Y) \In \gL(X_\cP,Y_\cP)$ for all $X$ and $Y$ in $\sC$.
    \end{description}

    \item ({\bf ``Bornologies''}) A family $\cM=\{\cM(Y): Y \in \sC\}$ of convex vector   bornologies associated to each space $Y$ in $\sC$ is called a {\em generating bornology} if
    \begin{description}
    \item[GB$_1$]  $\cM_\fin(Y) \In \cM(Y) \In \cM_\von(Y)$ for all $Y$ in $\sC$; and
    \item[GB$_2$]  $\gL(X,Y) \In L^\times(X^\cM,Y^\cM)$ for all $X$ and $Y$ in $\sC$.
    \end{description}
\end{enumerate}
\end{definition}

Classical examples of these notions are the ideals $\gK_p$ of precompact operators and $\gP$ of
absolutely summing operators (see \eg\cite{Pie80}), the generating systems $\cP_{pc}$ of precompact
topologies (see \eg\cite{Ran72}) and $\cP_{pn}$ of prenuclear topologies (see \eg\cite[p.\ 90]{Sch71}),
and the generating systems $\cM_{pc}$ of precompact bornologies and $\cM_{pn}$ of prenuclear bornologies (see \eg\cite{HN81}),
respectively. An interesting fact about these examples is that
we can visualize the notions of ``operators'', ``topologies'' and ``bornologies'' as vertices
of a triangle, and they can be transformed to each other by actions represented
as linking edges of the triangle.

\begin{definition}
Let $\gA$ be an operator ideal, $\cP$ a generating topology and $\cM$  a generating bornology on $\sC$.
\begin{enumerate}
    \item ({\bf ``Operators'' $\rightarrow$ ``Topologies''})
For each $X_0$ in $\sC$, the $\gA$--{\em topology} of $X_0$, denoted by $\sT(\gA)(X_0)$, is the projective topology of $X_0$ with respect to the family
$$
\{T \in \gA(X_0,Y):Y\in\sC\}.
$$
  In other words, a seminorm $p$ of $X_0$ is $\sT(\gA)(X_0)$--continuous if and only if there is a $T$ in $\gA(X_0,Y)$ for some $Y$ in
$\sC$ and a continuous seminorm $q$ of $Y$ such that
$$
p(x) \leq q(Tx), \quad\forall x \in X_0.
$$
In this case, we call $p$ an $\gA$--{\em seminorm} of $X_0$.

    \item ({\bf ``Operators'' $\rightarrow$ ``Bornologies''})
For each $Y_0$ in $\sC$, the $\gA$--{\em bornology} of $Y_0$, denoted by $\sB(\gA)(Y_0)$, is the inductive bornology of $Y_0$ with respect to the family
$$
\{T \in \gA(X,Y_0):X\in\sC\}.
$$
  In other words, a subset $B$ of $Y_0$ is $\sB(\gA)(Y_0)$--bounded if and only if there is a $T$ in $\gA(X,Y_0)$ for some $X$ in
$\sC$ and a topologically bounded subset $A$ of $X$ such that
$$
B \In TA.
$$
In this case, we call $B$ an $\gA$--{\em bounded subset} of $Y_0$.

    \item ({\bf ``Topologies'' $\rightarrow$ ``Operators''})
For $X$, $Y$ in $\sC$, let
$$
\sO(\cP)(X,Y) = \gL(X_\cP,Y)
$$
and
$$
\sO^b(\cP)(X,Y) = \gL^b(X_\cP,Y)
$$
be the vector space of all continuous operators from $X$ into $Y$ which is still continuous with
respect to the $\cP(X)$--topology, and which send a $\cP(X)$--neighborhood of zero to a bounded set, respectively.

    \item ({\bf ``Bornologies'' $\rightarrow$ ``Operators''})
For $X$, $Y$ in $\sC$, let
$$
\sO(\cM)(X,Y) = \gL(X,Y)\cap L^\times(X,Y^\cM)
$$
and
$$
\sO^b(\cM)(X,Y) = \gL^b(X,Y^\cM)
$$
be the vector space of all continuous operators from $X$ into $Y$ which send bounded
sets to $\cM(Y)$--bounded sets, and which send a neighborhood of zero to an $\cM(Y)$--bounded set, respectively.

    \item ({\bf ``Topologies'' $\leftrightarrow$ ``Bornologies''})
For $X$, $Y$ in $\sC$, the $\cP^\circ(Y)$--{\em bornology} of $Y$ (resp.\ $\cM^\circ(X)$--{\em topology} of $X$) is
defined to be the bornology (resp.\ topology) polar to
$\cP(X)$ (resp.\ $\cM(Y)$).  More precisely,
\begin{itemize}
    \item a bounded subset $A$ of $Y$ is $\cP^\circ(Y)$--bounded if and only if  its   polar
$A\po$ is a $\cP(Y'_\beta)$--neighborhood of zero; and
    \item a neighborhood $V$ of zero of $X$ is a $\cM^\circ(X)$--neighborhood of zero if and only if
$V^\circ$ is $\cM(X'_\beta)$--bounded.
\end{itemize}
\end{enumerate}
\end{definition}

\begin{theorem}\label{thm:tri_edges}
Let $\gA$ be an operator ideal, $\cP$ a generating topology and $\cM$ a generating bornology on $\sC$.  We have
\begin{enumerate}
    \item $\sT(\gA) = \{\sT(\gA)(X) : X \in \sC\}$ is a generating topology on $\sC$.
    \item $\sB(\gA) = \{\sB(\gA)(Y) : Y \in \sC\}$ is a generating bornology on $\sC$.
    \item $\sO(\cP) = \{\sO(\cP)(X,Y) : X,Y \in \sC\}$ is an operator ideal on $\sC$.
    \item $\sO^b(\cP) = \{\sO^b(\cP)(X,Y) : X,Y \in \sC\}$ is an operator ideal on $\sC$.
    \item $\sO(\cM) = \{\sO(\cM)(X,Y) : X,Y \in \sC\}$ is an operator ideal on $\sC$.
    \item $\sO^b(\cM) = \{\sO^b(\cM)(X,Y) : X,Y \in \sC\}$ is an operator ideal on $\sC$.
    \item $\cP^\circ = \{\cP^\circ(Y) : Y \in \sC\}$ is a generating bornology on $\sC$.
    \item $\cM^\circ = \{\cM^\circ(Y) : Y \in \sC\}$ is a generating topology on $\sC$.
\end{enumerate}
\end{theorem}
\begin{proof}
(1)--(6), together with the Banach space version of
(7) and (8), are done in \cite{WW88}.   For the locally convex space version of (7), we first note that (GB$_1$) follows
from (GT$_1$) and the bipolar theorem. To check (GB$_2$), let $X$ and $Y$ be
LCS's and $T\in \gL(X,Y)$.  Let $B$ be a $\cP^\circ(X)$--bounded subset of $X$ and we want to see that
$TB$ is $\cP^\circ(Y)$--bounded in $Y$.  Since $B^\circ$ is a $\cP(X'_\beta)$--neighborhood of zero of the
strong dual $X'_\beta$ of $X$, $(TB)^\circ = (T')^{-1}B^\circ$ is a $\cP(Y'_\beta)$--neighborhood of
zero of $Y'_\beta$ as a consequence of (GT$_2$) and the fact that $T'\in \gL(Y'_\beta,X'_\beta)$.  Hence, $TB$ is
$\cP^\circ(Y)$--bounded in $Y$, as asserted.

Finally, for (8) we note that (GT$_1$) is plain.  For (GT$_2$), let $X$ and $Y$ be LCS's and $T\in \gL(X,Y)$.  Let $V$ be a $\cM^\circ(Y)$-neighborhood of zero of $Y$ and we want to see that $T^{-1}V$
is an $\cM(X)$--neighborhood of zero of $X$.  Since $V^\circ$ is $\cM(Y'_\beta)$--bounded in $Y'_\beta$,
$(T^{-1}V)^\circ = T'V^\circ$ is a $\cM(X'_\beta)$--bounded subset of $X'_\beta$, as asserted.
\end{proof}

\begin{remark}
A seemingly more general setting  is to define for generating topologies $\cP$ and $\cP_1$, and generating
bornologies $\cM$ and $\cM_1$ the operator ideals with components $\sO(\cP/\cP_1)(X,Y) = \gL(X_\cP,Y_{\cP_1})$, $\sO(\cM/\cM_1)(X,Y)=L^\times(X^\cM,Y^{\cM_1})\cap\gL(X,Y)$ and
$\sO(\cP/\cM)(X,Y) = \gL(X_\cP,Y^\cM)$.  However, they will not give rise to new tools to us.     In fact,
we have $\sO(\cP/\cP_1) = \sO(\cP_1)^{-1}\circ\sO(\cP)$ \cite{Step80}, $\sO(\cM/\cM_1) = \sO(\cM)\circ\sO(\cM_1)^{-1}$ \cite{Step83}, and $\sO(\cP/\cM) = \sO^b(\cM)\circ\sO(\cP) = \sO(\cM)\circ\sO^b(\cP)$.   Readers are referred to Pietsch's classic \cite{Pie80} for information regarding quotients and products of operator ideals.
\end{remark}

Let $p$ be a continuous seminorm
of a LCS $X$ and $B$ an absolutely convex bounded subset of a LCS $Y$.  Denote by
$X_p$ the normed space $X/p^{-1}(0)$ equipped with norm $\|x + p^{-1}(0)\| = p(x)$, and
by $Y(B)$ the normed space $\bigcup_{\lambda>0} \lambda B$ equipped with norm $r_B(x) =
\inf\{\lambda > 0 : x \in \lambda B\}$.  Let $\wX_p$ be the completion of $X_p$.  Define $Q_p :
X \lrw X_p$, $\wQ_p : X \lrw \wX_p$ and $J_B : Y(B) \lrw Y$ to be the canonical maps.

\begin{theorem}[\cite{W94}]\label{T:TB}
Let $\gA$ be an operator ideal on LCS's.    We have
\begin{enumerate}
    \item A continuous seminorm $p$ of $X$ is an $\gA$--seminorm if and only if $Q_p \in \gA^\inj(X,X_p)$ if and only if
$\wQ_p \in \gA^\inj(X,\wX_p)$.
    \item A bounded disk $B$ of $Y$ is an $\gA$--bounded set if and only if $J_B \in \gA^\bsur(Y(B),Y)$.  Whenever $\gA$
is surjective, we can replace $\gA^\bsur$ by $\gA^\sur$.
\end{enumerate}
\end{theorem}

For operator ideals $\gA$ on Banach spaces, Stephani \cite{Step70, Step73} achieved that
$\sO(\sT(\gA)) = \gA^\inj$ and $\sO(\sB(\gA)) = \gA^\sur$.  However,
we have two constructions $\sO$ and $\sO^b$ in the context of LCS's.  Unlike the Banach space version, they give rise to different ideals.  For example,
let $\cM_{pc}$ be the generating system of precompact bornologies (\ie the bornologies determined by
totally bounded convex sets).  Then $\gK_p = \sO^b(\cM_{pc})$ is the ideal of precompact operators (\ie
those sending a neighborhood of zero to a totally bounded set) and $\gK_p^{loc} = \sO(\cM_{pc})$ is the
ideal of locally precompact operators (\ie those sending bounded sets to totally bounded sets).  Randtke \cite{Ran72} indicated that $\gK_p(X,Y) = \gK_p^{loc}(X,Y)$ holds for all LCS $Y$ if and only if $X$ is a Schwartz space.  On the other hand, it is straightforward to make the following observation.

\begin{proposition}\label{P:TOTBOB}
For a generating topology $\cP$ and a generating bornology $\cM$ on LCS's,
$\sO(\cP)$ and $\sO^b(\cP)$ give
rise to the same ideal  topology, namely
$$
\sT(\sO(\cP)) = \sT(\sO^b(\cP)) = \cP,
$$
and $\sO(\cM)$ and $\sO^b(\cM)$ give rise to the same ideal bornology, namely
$$
\sB(\sO(\cM)) = \sB(\sO^b(\cM)) = \cM.
$$
Moreover, $\sO(\cP)$ and $\sO^b(\cP)$ are injective, $\sO(\cM)$ is
bornologically surjective and $\sO^b(\cM)$ is surjective.
\end{proposition}

\begin{proposition}\label{P:OTOOBO}
Let $\gA$ be an operator ideal on LCS's.  We have
\begin{enumerate}
    \item $\sO^b(\sT(\gA)) \In \gA^\inj \In \sO(\sT(\gA))$.
    \item $\sO^b(\sB(\gA)) \In \gA^\bsur \In \sO(\sB(\gA))$.
\end{enumerate}
\end{proposition}
\begin{proof}
Let $T$ be a (topologically) bounded linear operator from a LCS $X$ into a LCS $Y$,
\ie $T\in \gL^b(X,Y)$.  Then there is a continuous seminorm $p$ of $X$ and an absolutely convex bounded
subset $B$ of $Y$ such that $T$ sends $V_p = \{x\in X: p(x) \leq 1\}$ into $B$.  It is plain that
$T$ has a decomposition
$$
\begin{CD}
X@>T>>Y\\
@VQ_pVV @AAJ_BA\\
X_p@>>T_0> Y(B),
\end{CD}
$$
where $T_0 \in \gL(X_p,Y(B))$ is the unique bounded operator induced by $T$.

If $T\in\sO^b(\sT(\gA))(X,Y)$ then $p$ can be chosen to be an $\gA$--seminorm of $X$.  By Theorem
\ref{T:TB}, $Q_p\in \gA^\inj(X,X_p)$ and hence $T=J_BT_0Q_p \in \gA^\inj(X,Y)$.    Similarly, if
$T\in \sO^b(\cM))(X,Y)$ then $B$ can be chosen to be an $\gA$--bounded subset of $Y$.   By
Theorem \ref{T:TB} again, $J_B\in \gA^\bsur(Y(B),Y)$ and hence $T=J_BT_0Q_p\in \gA^\bsur(X,Y)$.
In other words, $\sO^b(\sT(\gA)) \In \gA^\inj$ and $\sO^b(\sB(\gA)) \In \gA^\bsur$.  The
other inclusions follows from the injectivity of $\sO(\sT(\gA))$ and the bornological surjectivity
of $\sO(\sT(\gA))$.
\end{proof}

\begin{proposition}\label{P:T2B}
Let $\cP$ be a generating topology on LCS's.  If the operator ideal $\gA = \sO(\cP)$
is symmetric (resp.\ $\gA = \sO^b(\cP)$ is symmetric) then
$$
\cP^\circ(Y) = \sB(\sO(\cP))(Y) \quad (\mbox{\resp} \cP^\circ(Y) = \sB(\sO^b(\cP))(Y))
$$
for all infrabarrelled LCS $Y$.
\end{proposition}
\begin{proof}
Let $B$ be a bounded disk in $Y$.  Suppose firstly that
$B$ is  $\gA$-bounded.  Then there is a normed space $N$ such
that $TU_N \Sset B$.  Hence $B^\circ \Sset (T')^{-1}U_{N'}$.  Now, the symmetry of $\gA$ implies
$T' \in \gA(Y'_\beta,N')$.  Thus, $B^\circ$ is an $\gA$-neighborhood of zero of $Y'_\beta$.
It follows from $\sT(\gA)
= \sT(\sO(\cP)) = \sT(\sO^b(\cP)) = \cP$ that $B^\circ$ is $\cP(Y'_\beta)$-neighborhood of zero of $Y'_\beta$.  Hence
$B$ is $\cP^\circ(Y)$-bounded.

Conversely,  assume that $B$ is $\cP^\circ$-bounded in $Y$.    In other words,
$B^\circ$ is an $\gA$-neighborhood of zero of $Y'_\beta$.  Therefore, there is a Banach space $F$ and
a $T$ in $\gA(Y'_\beta,F)$ such that $B^\circ \Sset T^{-1}U_F$.  Hence the second polar $B^{\circ\bullet}$  of
$B$ in $Y''_{\beta\beta}$ is $\gA$-bounded since $B^{\circ\bullet} \in T'U_{F'}$ and $T' \in \gA(F',Y''_{\beta\beta})$.  Let $K_Y$ be the canonical embedding of $Y$ into $Y''_{\beta\beta}$.
The infrabarrelledness of $Y$ ensures that $K_Y$ is a topological injection.  As a result,
the inclusion $K_YB \In B^{\circ\bullet}$ establishes the existence of a $k_B$ in $\gL(Y(B),Y''_{\beta\beta}(B^{\circ\bullet}))$ such that
$J_{B^{\circ\bullet}}k_B = K_XJ_B$.  Then $J_B \in (\gA^\bsur)^\inj(Y(B),Y)$ because $J_{B^{\circ\bullet}}\in
\gA^\bsur(Y''_{\beta\beta}({B^{\circ\bullet}}), Y''_{\beta\beta})$ by Theorem \ref{T:TB}.  However,
$(\gA^\bsur)^\inj = (\gA^\inj)^\bsur = \gA^\bsur$ since $\gA = \sO(\cP)$ (or $\gA = \sO^b(\cP)$) is always injective.
This implies that $B$ is $\gA$-bounded, \ie $\sB(\sO(\cP))$-bounded in $Y$, by Theorem \ref{T:TB} again.
\end{proof}

\begin{proposition}\label{P:B2T}
Let $\cM$ be a generating bornology on LCS's.  If the operator ideal $\gA = \sO(\cM)$
is symmetric (resp.\ $\gA = \sO^b(\cM)$ is symmetric) then
$$
\cM^\circ(X) = \sT(\sO(\cM))(X) \quad (\mbox{\resp} \cM^\circ(X) =\sT(\sO^b(\cM))(X))
$$
for all infrabarrelled LCS $X$.
\end{proposition}
\begin{proof}
Let $V$ be a  closed, absolutely convex neighborhood of zero of $X$.
Suppose firstly that $V$ is an $\gA$-neighborhood of $X$ then there is a normed space
$N$ and an $T$ in $\gA(X,N)$ such that $T^{-1}U_N \In V$.  Hence $V^\circ \In T'U_{N'}$ and
thus $V^\circ$ is $\cM$-bounded in $X'_\beta$ since $T' \in \gA(N',X'_\beta)$.  So $V$ is an
$\cM^\circ(X)$-bounded subset of $X$.

Conversely, assume that $V$ is an $\cM^\circ$-neighborhood of zero of $X$.   Then $V^\circ$ is
$\cM(X'_\beta)$-bounded in the strong dual space $X'_\beta$ of $X$.  Hence there is a normed space
$N$ and an $T$ in $\gA(N,X'_\beta)$ such that $TU_N \Sset V^\circ$.  Consequently, ${V^{\circ\bullet}}
\In (T')^{-1}U_{N'}$ and thus ${V^{\circ\bullet}}$ is an $\gA$-neighborhood of zero of $X''_{\beta\beta}$
as $T' \in \gA(X''_{\beta\beta}, N')$.  Since $X$ is infrabarrelled, $K_X$ is continuous.
By (GT$_2$), $V = K_X^{-1}  {V^{\circ\bullet}} = {V^{\circ\bullet}}\cap X$ is an $\gA$-neighborhood
of zero of $X$.
\end{proof}

\definition
A generating topology $\cP$ on $LCS$'s is said to
have the
{\it subspace property\/} if whenever $Y$ is a subspace of a $LCS\ X$,
$Y_{\cP}$ is also a subspace of $X_{\cP}$, \ie the $\cP$--topology of
$Y$
coincides with the subspace topology inherited from the $\cP$--topology
of $X$.
See Jarchow \cite{Jar81} for the
Banach space version.

Let $\gA$ be an operator ideal on LCS's or Banach spaces.
$\gA^{\dual}$
denotes the operator ideal with components
$$
\gA^{\dual}(X,Y)=\{T\in\sL(X,Y) :T'\in\gA(Y'_\b,X'_\b)\}.
$$

\begin{proposition}\label{P:TBdual}
Let $\cP$ be a generating topology on $LCS$'s
and $X$
be an infrabarrelled $LCS$. Then
\begin{enumerate}[(a)]
    \item  $\sO^b(\cP)^{\dual}(X,Y)=\sO^b(\cP^{\circ})(X,Y),\qquad\forall
LCS\ Y$.
    \item   $\sO^b(\cP^\circ)^{\dual}(X,Y)\subseteq\sO^b(\cP)(X,Y),\qquad
\forall LCS\ Y$.
\end{enumerate}
If, in addition, $\sO(\cP)$ is symmetric or $\cP$ has the subspace
property then
\begin{enumerate}[{(b)$'$}]
\item
$\sO^b(\cP^{\circ})^{\dual}(X,Y)=\sO^b(\cP)(X,Y),\qquad
\forall LCS\ Y$.
\end{enumerate}
\end{proposition}
\begin{proof} (a) Let $T\in\sO^b(\cP)^{\dual}(X,Y)$, \ie
$T'\in\sO^b(\cP)
(Y'_\b,X'_\b)$. Then there is a $\cP(Y'_\b)$--neighbor\-hood $V$ of $0$
in
$Y'_\b$ such that $T'V$ is bounded in $X'_\b$. Hence $U=(T'V)^\circ=
T^{-1}V^\circ$ is a closed bornivorous barrel in $X$. Since $X$ is
infrabarrelled, $U$
is a  $0$--neighborhood in $X$. Now $TU\subseteq V^\circ$ ensures that
$T\in\sO^b
(\cP^\circ)(X,Y)$. Conversely, if $T\in\sO^b(\cP^\circ)(X,Y)$ then
there is a
$0$--neighborhood $U$ in $X$ such that $A=TU$ is
$\cP^\circ(Y)$--bounded in
$Y$. Hence $A^\circ=(T')^{-1}U^\circ$ is a $\cP(Y'_\b)$--neighborhood
of $0$
in $Y'_\b$. Now $T'A^\circ\subseteq U^\circ$ implies that
$T'\in\sO^b(\cP)
(Y'_\b,X'_\b)$.

(b) Let $T\in\sO^b(\cP^\circ)^{\dual}(X,Y)$, \ie
$T'\in\sO^b(\cP^\circ)(Y'_\b,
X'_\b)$. Then there is a $0$--neighborhood $V$ in $Y'_\b$ such that
$T'V$ is
$\cP^\circ(X'_\b)$--bounded in $X'_\b$. Hence $U=(T'V)^\bullet$ is a
$\cP(X''_{\b\b})$--neighborhood of $0$ in $X''_{\b\b}$. Let
$U_0=K^{-1}_XU$.
By the functorial property $(GT_2)$ of $\cP$, $U_0$ is a $\cP(X)$--neighborhood
of $0$ in $X$. It is easy to see that $TU_0\subseteq V^\circ$ and thus
$T\in\sO^b
(\cP)(X,Y)$.

(b)$'$ Assume, in addition to those in (b), that $\gA=\sO(\cP)$ is
symmetric.
Let $T\in\sO^b(\cP)(X,Y)$. We want to verify that
$T'\in\sO^b(\cP^\circ)
(Y'_\b,X'_\b)$. By assumption, there is a $\cP(X)$--neighborhood $U$ of
$0$ in
$X$ such that $A=TU$ is bounded in $Y$. Now $A^\circ=(T')^{-1}U^\circ$
suggests us to check if $U^\circ$ is $\cP^\circ(X'_\b)$--bounded in
$X'_\b$. Since
$\cP=\sT(\sO(\cP))=\sT(\gA)$, there is a Banach space $F$ and an $R$ in
$\gA(X,F)$
such that $U\supseteq R^{-1}U_{F}$. Therefore, $U^\circ\subseteq
R'U_{F'}$ and thus
$U^{\circ\bullet}\supseteq((R')')^{-1}U_{F''}$ where $(R')'$ is the
double
adjoint of $R$ from $X''_{\b\b}$ into $F''$. Since $\gA$ is symmetric,
$(R')'\in\gA(X''_{\b\b},F'')$ and thus
$U^{\circ\bullet}$ is a $\cP(X''_{\b\b})$--neighborhood of $0$ in
$X''_{\b\b}$. Consequently, $U^\circ$ is a $\cP^\circ(X'_\b)$--bounded
subset
of $X'_\b$, as asserted.

Finally, if the subspace property of $\cP$ is assumed instead of the
symmetry
of $\sO(\cP)$ then the $\cP(X)$--neighborhood $U$ of $0$ in $X$ above
is
induced from a $\cP(X''_{\b\b})$--neighborhood $V$ of $0$ in
$X''_{\b\b}$,
\ie $K_XU=V\cap K_XX$ and thus $U^\circ=V^\bullet$ is $\cP^\circ
(X'_\b)$--bounded in $X'_\b$, where $K_X$ is the evaluation map from
$X$ into
$X''_{\b\b}$.
\end{proof}


\section{LCS's defined by operators, topologies and bornologies}\label{s:Groth-spaces}

\begin{theorem}\label{thm:3.4} Let $\gA$ be an operator ideal on LCS's and
$X$ be a
LCS. The following are all equivalent.
\begin{enumerate}[(1)]
\item $X$ is $\gA$--topological.
\item For each continuous seminorm $p$ on $X$, $Q_p\in\gA^{\inj}(X,X_p)$, or equivalently, $\wQ_p\in\gA^{\inj}(X,\wX_p)$.
\item $\sL^b(X,Y)\subseteq\gA^{\inj}(X,Y)$ for every LCS $Y$.
\item $\sL(X,F)=\gA^{\inj}(X,F)$ for every normed (or Banach)
space $F$.
\item $id_X\in\sL(X_{\gA},X)$, where $X_{\gA}$ is the LCS $X$
equipped
with the $\gA$--topology.
\end{enumerate}
\end{theorem}
\begin{proof} (1)$\Leftrightarrow$(2) is contained in Theorem \ref{T:TB}.
(1)$\Leftrightarrow$(5) and (2)$\Rightarrow$(3)$\Rightarrow$(4) are
trivial.
(4)$\Rightarrow$(1) is due to Theorem \ref{T:TB} again.
\end{proof}

In the following, $\sC$ denotes either the class of all LCS's or the
class of
all Banach spaces. The next result is a generalization of a result of
Jarchow
\cite[Proposition 3]{Jar84}.

\begin{theorem}[\cite{W94}]\label{thm:3.5} Let $\gA$ be a surjective operator ideal on
$\sC$.
If $X$, $Y\in\sC$ and $Y$ is a (topological) quotient space of $X$ then
the
$\gA$--topology of $Y$ is the quotient topology induced by the
$\gA$--topology
of $X$.
In particular,
a quotient space of an $\gA$--topological space is again an
$\gA$--topological space.
\end{theorem}


\begin{theorem}\label{thm:3.10} Let $\gA$ be an operator ideal on LCS's and
$Y$ be a
LCS. The following are all equivalent.
\begin{enumerate}[(1)]
\item $Y$ is $\gA$--bornological.
\item $J_B\in\gA^{\bsur}(Y(B),Y)$ for each bounded disk $B$ in
$Y$.
\item $\sL^b(X,Y)\subseteq\gA^{\bsur}(X,Y)$ for every  LCS $X$.
\item $\sL(N,Y)=\gA^{\bsur}(N,Y)$ for every normed space $N$.
\item $id_Y\in L^{\times}(Y,Y^{\gA})$, where $Y^{\gA}$ is the
convex
bornological vector space $Y$ equipped with the $\gA$--bornology.
\end{enumerate}
\noindent
In case $Y$ is infracomplete they are all equivalent to
\begin{enumerate}[(4)$'$]
\item $\sL(E,Y)=\gA^{\bsur}(E,Y)$ for every Banach space $E$.
\end{enumerate}
If $\gA$ is surjective we can replace $\gA^{\bsur}$ by $\gA$ in all of
the
above statements.
\end{theorem}
\begin{proof} (1)$\Leftrightarrow$(5) is by definition. It is plain that
(2)$\Rightarrow$(3)$\Rightarrow$(4)$\Rightarrow$(4)$'$.
(4)$\Rightarrow$(1) (or
(4)$'\Rightarrow1$ in case $Y$ is infracomplete) and (1)$\Leftrightarrow$(2) is due to Theorem \ref{T:TB}. The last
assertion is
a consequence of Proposition \ref{hulls}.
\end{proof}

Let $\sC$ be either the class of all LCS's or the class of all Banach
spaces.

\begin{theorem}\cite{W94}\label{thm:3.11} Let $\gA$ be an injective operator ideal on
$\sC$ and $X, Y \in \sC$.
If $Y$ is a (topological) subspace of $X$ then the
$\gA$--bornology of $Y$ is the subspace bornology inherited from the
$\gA$--bornology of $X$.
In particular, a subspace of an $\gA$--bornological space is again an
$\gA$--bornological
space.
\end{theorem}


\begin{theorem}\label{thm:6.4} Let $\gA$ be an operator ideal on LCS's. Let
$\sG=\sT(\gA)$ be the ideal topology on LCS's generated by $\gA$. A
LCS
$X$ is $\gA$--topological if and only if
$$
\sL^b(X,Y)\cap\sO(\sG)(X,Y)=\sO^b(\sG)(X,Y)
$$
for each LCS $Y$.
\end{theorem}
\begin{proof} By Theorem \ref{thm:3.4}, if $X$ is $\gA$--topological then
$\sO^b(\sG)(X,Y)
=\sL^b(X,Y)$ and $\sO(\sG)(X,Y)=\sL(X,Y)$. The equality follows.
Conversely,
assume the equality holds for every LCS $Y$. It suffices to show that
$\sL
(X,N)\subseteq\sO^b(\sG)(X,N)$ for each normed space. Let $N_{\sG}$ be
the
LCS given by equipping $N$ with the $\sG(N)$--topology. By the
functorial
property $(GT_2)$ of $\sG$, any $T$ in $\sL(X,N)$ also belongs to
$\sL(X_{\sG},
N_{\sG})=\sO(\sG)(X,N_{\sG})$. By the hypothesis,
$T\in\sO^b(\sG)(X,N_{\sG})$.
Since $\sG(N)$ is compatible with the dual pair $(N,N')$ by $(GT_1)$,
we have
$T\in\sO^b(\sG)(X,N)$. It follows the desired
assertion.
\end{proof}

\begin{remark} If we let $\cM=\sB(\gA)$ then a LCS $Y$ being $\gA$--bornological implies
$$
\sL^b(X,Y)\cap\sO(\cM)(X,Y)=\sO^b(\cM)(X,Y)
$$
for each LCS $X$.  We do not know if the converse is true.
\end{remark}

Let $\gA$ be an operator ideal on LCS's.  Denote by
$\gA_{\bb}$ the operator ideal defined on Banach spaces such that
$\gA_{\bb}(E,F)=\gA(E,F)$ for every pair $E$ and $F$ of Banach spaces.
Conversely, let $\gA$ be an operator ideal on Banach spaces. There are
many ways to extend $\gA$ to an operator ideal
$\gA_0$ on LCS's
in the sense that $(\gA_0)_{\bb} = \gA$.
In \cite{Pie80},
Pietsch mentioned six different ways to extend $\gA$ to an operator
ideal on
LCS's. Among them, we are interested in
\begin{eqnarray*}
\gA^{\inf} &=& \{RS_0 T\in \gL(X,Y): T\in\gL(X,X_0), S_0\in \gA(X_0,Y_0), R\in\gL(Y_0,Y)\},\\
\gA^{\rup} &=&  \{S\in \gL(X,Y):  \forall B\in\gL(Y,Y_0), \exists  A\in\gL(X,X_0),  S_0\in\gA(X_0,Y_0)
    \text{ such that } BS=S_0 A\},\\
\gA^{\lup} &=& \{S\in \gL(X,Y):   \forall B\in\gL(X_0,X), \exists A\in\gL(Y_0,Y),  S_0\in\gA(X_0,Y_0)
    \text{ such that } SB=AS_0\},\\
\gA^{\sup} &=& \{S\in \gL(X,Y): RST\in \gA(X_0,Y_0), \text{ for all } T\in\gL(X_0,X) \text{ and } R\in\gL(Y,Y_0)\}.
\end{eqnarray*}
Here, $X,Y$ run through all LCS's and $X_0,Y_0$ run through all Banach spaces.

\begin{definition}[\cite{W94}]
 Let $\gA$ be an operator ideal on Banach spaces. We
call a
continuous seminorm $p$ on a LCS $X$ a $\Groth(\gA)$--{\it seminorm\/}
if
there is a continuous seminorm $q$ on $X$ such that $p\leq q$ and
$\wQ_{pq}
\in\gA(\wX_q,\wX_p)$.
The $\Groth(\gA)$--{\it topology\/} on $X$ is defined to be the
locally
convex (Hausdorff) topology on $X$ which has a subbase determined by
all $\Groth(\gA)$--seminorms.

A LCS $X$ is  a $\Groth(\gA)$--space if its topology coincides with the
$\Groth(\gA)$--topology.  It is equivalent to say that the identity map
$id_X\in\gA^\rup(X,X)$.
\end{definition}

Let $\sG$ be a generating topology on Banach spaces. Define
$\sG^{\bl}(X)$ on
each LCS $X$ to be the coarsest locally convex (Hausdorff) topology on
$X$
among those $\sG_0(X)$ such that the inclusion
$$
\sL(X,F)\subseteq\sL(X_{\sG_0},F_{\sG})
$$
holds for every Banach space $F$.
It is clear that for each LCS $X$,
$\sG_\s(X)\leq\sG^{\bl}(X)\leq\sG_{\ori}
(X)$ and a continuous seminorm $p$ on $X$ is $\sG^{\bl}(X)$--continuous
if and
only if there is a Banach space $F$, an $S$ in $\sL(X,F)$ and a
$\sG(F)$--continuous seminorm $r$ on $F$ such that $p(x)\leq r(Sx)$ for
all
$x$ in $X$.

\begin{lemma}\label{lem:4.8} $\sG^{\bl}=\{\sG^{\bl}(X) :X\ LCS\}$   is the
minimal extension of $\sG$ to LCS's.
\end{lemma}
\begin{proof} It is easy to see that $\sG^{\bl}$ is a generating
topology on
LCS's. Let $E$ be a Banach space. By definition of $\sG^{\bl}$,
$\sG^{\bl}
(E)\leq\sG(E)$. On the other hand,
$id_E\in\sL(E,E)\subseteq\sL(E_{\sG^{\bl}},
E_{\sG})$ implies $\sG^{\bl}(E)\geq\sG(E)$. So $\sG^{\bl}$ is an
extension of
$\sG$ to LCS's. The minimality of $\sG^{\bl}$ is
obvious.
\end{proof}

\begin{theorem}\label{thm:4.9} Let $\gA$ be an operator ideal on Banach spaces.
The
minimal extension $\sG^{\bl}$ of $\sG=\sT(\gA)$ coincides with the
$\Groth(\gA)$--topology.
\end{theorem}
\begin{proof} Without loss of generality, we can assume that $\gA$ is injective since
$\sT(\gA)=\sT(\gA^{\inj})$ by Theorem \ref{T:TB}.
Let $p=r\circ S$
be a $\sG^{\bl}$--continuous seminorm on a LCS $X$ where
$S\in\sL(X,F)$ and $r$ is a $\sG$--continuous seminorm on a Banach
space $F$. Then we have $\V\wQ_p(x)\V_{\wX_p}=\V\wQ_r(Sx)\V_{\wF_r}$
for all $x$ in $X$. It follows
that there is an isometry $S_0$ in $\sL(\wX_p,\wF_r)$ such that
$S_0\wQ_p=\wQ_rS$.
Note that $\wQ_r\in\gA(F,\wF_r)$. Define a continuous seminorm $q$ on
$X$ by
$q(x)=\V\wQ_r\V\ \V Sx\V$.
Now $q(x)\geq\V\wQ_rSx\V=p(x)$ and
we have an $S_2$ in $\sL(\wX_q,F)$ induced by $S$.
Since $S_0$ is an injection and
$S_0\wQ_{pq}=\wQ_rS_2\in\gA(\wX_q,
\wF_r)$, $\wQ_{pq}\in\gA(\wX_q,\wX_p)$, \ie $p$ is a
$\Groth(\gA)$--seminorm.

Conversely, if $p$ is a $\Groth(\gA)$--continuous seminorm on $X$ then
there
is a continuous seminorm $q$ on $X$ with $p\leq q$ such that $\wQ_{pq}
\in\gA(\wX_q,\wX_p)=\sL((\wX_q)_{\sG},\wX_p)$ by Theorem \ref{thm:Bspaceversion}(4a). In
other words, the seminorm $r$ on $\wX_q$ defined by $r(y)=\V\wQ_{pq}(y)
\V_{\wX_p}$, $y\in\wX_q$, is $\sG$--continuous. Note that
$\wQ_p=\wQ_{pq}\wQ_q$
implies that $p(x)=\V\wQ_{pq}\wQ_q(x)\V=r(\wQ_qx)$. It simply says that
$p$ is
a $\sG^{\bl}$--continuous seminorm.
\end{proof}

\begin{theorem}\label{thm:4.14} Let $\gA$ be an operator ideal on Banach spaces
with $\sG=\sT(\gA)$. Then $\sO(\sG^{\bl})=(\gA^{\inj})^{\rup}$.
\end{theorem}
\begin{proof} Let $X$ and $Y$ be LCS's.
Assume $T\in\sO
(\sG^{\bl})(X,Y)$. Then for every Banach space $F$ and $S$ in
$\sL(Y,F)$, $ST\in\sO
(\sG^{\bl})(X,F)=\sL(X_{\sG^{\bl}},F)$. Hence there is a
$\sG^{\bl}$--continuous seminorm $p$ on $X$ such that $\V STx\V\leq
p(x)$. By
Theorem \ref{thm:4.9}, there is a continuous seminorm $q$ on $X$ such that $p\leq
q$
and $\wQ_{pq}\in\gA^{\inj}(\wX_q,\wX_p)$.
Let $R\in\gL(\wX_p,F)$ is induced by the inequality $\V STx\V\leq p(x)$.
It is then not difficult to see that $ST=R\wQ_{pq}\wQ_{q}$, and thus
 $T\in(\gA^{\inj})^{\rup}(X,Y)$.

Conversely, assume $T\in(\gA^{\inj})^{\rup}(X,Y)$. Then for every
continuous
seminorm $p$ on $Y$ there exists a Banach space $E$, an $R$ in
$\sL(X,E)$ and an
$S$ in $\gA^{\inj}(E,\wY_p)$ such that $\wQ_pT=SR$. Now
$S\in\sL(E_{\sG},\wY_p)$
and $R\in\sL(X,E)\subseteq\sL(X_{\sG^{\bl}},E_{\sG})$ imply $\wQ_pT=SR\in
\sL(X_{\sG^{\bl}},\wY_p)$. Since it is true for every continuous
seminorm $p$
on $Y$, $T\in\sL(X_{\sG^{\bl}},Y)$, \ie $T\in\sO(\sG^{\bl})
(X,Y)$.
\end{proof}

\begin{definition}[\cite{W94}] Let $\gA$ be an operator ideal on Banach spaces. A
bounded
$\s$--disk $A$ in a LCS $X$ is said to be $\Groth(\gA)$--{\it
bounded\/} in
$X$ if there is a bounded $\s$--disk $B$ in $X$ such that $A\subseteq B$
and the
canonical map $J_{BA}\in\gA(X(A),X(B))$. Note that, in this case, both
$X(A)$
and $X(B)$ are Banach spaces.
The
$\Groth(\gA)$--{\it bornology\/} on a LCS $X$ is defined to be the convex
vector
bornology
on $X$ with a subbase consisting of $\Groth(\gA)$--bounded $\s$--disks
in
$X$.

A LCS is a  {\it co}--$\Groth(\gA)$--{\it space,\/} if all bounded
$\s$--disks
in $X$ are $\Groth(\gA)$--bounded. It is equivalent to say that
$id_X\in
\gA^{\lup} (X,X)$.
\end{definition}

Let $\cM$ be a generating bornology on Banach spaces. We define, for
each
LCS $X$, a convex vector bornology $\cM^{\bl}(X)$ on $X$ to be the
smallest
convex (separated) vector bornology among those $\cM_0$ on $X$ such
that
$$
\sL(E,X)\subseteq L^{\times}(E^{\cM},X^{\cM_0})
$$
holds for every Banach space $E$.
It is easy to see that $\cM_{\fin}(X)\subseteq
\cM^{\bl}(X)\subseteq\cM_{\von}(X)$ and the family of subsets $B$ in $X$
in the
form of $B=TA$ for some $T\in\sL(E,X)$ and $\cM$--bounded set $A$ in a
Banach
space $E$ forms a basis of the bornology $\cM^{\bl}(X)$ for each LCS $X$.

\begin{lemma}\label{lem:4.12} Let $\cM$ be a generating bornology on Banach
spaces.
$\cM^{\bl}$ is the minimal extension of $\cM$ to LCS's.
\end{lemma}
\begin{proof} Similar to Lemma \ref{lem:4.8}.
\end{proof}

\begin{theorem}\label{thm:4.13} Let $\gA$ be an operator ideal on Banach
spaces. The
minimal extension $\cM^{\bl}$ of $\cM=\sB(\gA)$ coincides with the
$\Groth(\gA)$--bornology.
\end{theorem}
\begin{proof}
Without loss of generality, we can assume that $\gA$ is surjective since $\sB(\gA)=
\sB(\gA^{\sur})$ by Theorem \ref{T:TB}. Let $A$ be a $\Groth(\gA)$--bounded
$\s$--disk in a LCS $X$. By definition, there is a
bounded
$\s$--disk $B$ in $X$ such that $A\subseteq B$ and
$J_{BA}\in\gA(X(A),X(B))
=L^{\times}(X(A),X(B)^{\cM})$. In other words, $C=J_{BA}U_{X(A)}$ is
$\cM$--bounded in $X(B)$. Now $A\subseteq
J_AU_{X(A)}=J_BJ_{BA}U_{X(A)}=J_BC$ implies that
$A$ is $M^{\bl}$--bounded.

Conversely, if $A=SB$ is $\cM^{\bl}$--bounded in
$X$ with some $S$ in $\sL(E,X)$ and $\cM$--bounded $\s$--disk $B$ in a
Banach
space $E$. Let $C=\lam SU_E$ for some $\lam>0$ such that
$\lam U_E\supseteq
B$. We have $C\supseteq A$.
Let $S_0\in\gL(E(B),X(A))$ and $S_2\in\gL(E, X(C))$ be induced by $S$.
Since $B$ is
$\cM$--bounded in $E$, $J_B\in\gA(E(B),E)$ and $J_{CA}S_0=S_2J_B\in
\gA(E(B),X(C))$. Finally the surjectivity of $S_0$ ensures that
$J_{CA}\in\gA
(X(A),X(C))$.
\end{proof}

\begin{theorem}\label{thm:4.15} Let $\gA$ be an operator ideal on Banach spaces
with
$\cM=\sB(\gA)$. Then
$$
\sO(\cM^{\bl})(X,Y)\subseteq(\gA^{\sur})^{\lup}(X,Y),\quad\forall
LCS\text{'s }
X,Y.
$$
If $X$ is infracomplete (in particular, a Banach space) then we have
$$
\sO(\cM^{\bl})(X,Y)=(\gA^{\sur})^{\lup}(X,Y),\quad\forall LCS\ Y.
$$
\end{theorem}
\begin{proof} Similar to a previous theorem except that we shall use
Theorem
\ref{thm:4.13} instead of Theorem \ref{thm:4.9}. The introduction of the infracompleteness
is merely to
give us a chance to utilize the extension condition.
\end{proof}


We provide a new proof for the following result.

\begin{theorem}[\cite{W94}]\label{thm:4.16}
Let $\gA$ be an operator ideal on Banach
spaces. The Groth$(\gA^\inj)$--topology coincides with the
$\gA^{\rup}$--topology on every LCS, and the Groth$(\gA^\sur)$--bornology
coincides with the $\gA^{\lup}$--bornology on every infracomplete
LCS.  In particular, we have
\begin{enumerate}[(a)]
\item
A LCS $X$ is a Groth$(\gA^{\inj})$--space if and only if
$X$ is an
$\gA^{\rup}$--topological space.
\item
An infracomplete LCS $X$ is a
co--Groth$(\gA^{{{\sur}}})$--space if and only if $X$ is
an $\gA^{\lup}$--bor\-no\-logical space.
\item
The $\gA$--topology (resp.\ $\gA$--bornology) coincides with
the Groth$(\gA^\inj)$--topology (resp.\ Groth$(\gA^\sur)$--bornology)
on Banach spaces.
\end{enumerate}
\end{theorem}
\begin{proof} Let $\sG=\sT(\gA)$ and $\sM=\sB(\gA)$ be the ideal
topology and
the ideal bornology on Banach spaces generated by $\gA$, respectively. Let $p$ be a
continuous seminorm on a LCS $X$.  We observe the following
equivalences:
\begin{enumerate}
\item[{\mbox{ \quad}}] $p$ is a $\Groth(\gA^{\inj})$--continuous seminorm on $X$.
\item[{$\Leftrightarrow$}] $p$ is an $\sO(\sG^{\bl})$--continuous
seminorm on $X$ by Theorem \ref{thm:4.9}.
\item[{$\Leftrightarrow$}] $p$ is an $(\gA^{\inj})^{\rup}$--continuous
seminorm on $X$ by
Theorem \ref{thm:4.14}.
\item[{$\Leftrightarrow$}]
$\wQ_p\in[(\gA^{\inj})^{\rup}]^{\inj}(X,\wX_p)$ by Theorem \ref{T:TB}.
\item[{$\Leftrightarrow$}] $\wQ_p\in(\gA^{\rup})^{\inj}(X,\wX_p)$ by
\cite[Proposition  3.5]{WW93}.
\item[{$\Leftrightarrow$}] $p$ is an $\gA^{\rup}$--continuous seminorm
on $X$ by Theorem \ref{T:TB}.
\end{enumerate}

For the bornological case, assuming that $X$ is infracomplete, we have
for each bounded $\s$--disk $A$ in $X$:
\begin{enumerate}
\item[{\mbox{\ \quad}}] $A$ is a co--$\Groth(\gA^{\sur})$--bounded set in $X$.
\item[{$\Leftrightarrow$}] $A$ is an $\sO(\sM^{\bl})$--bounded set in
$X$ by Theorem
\ref{thm:4.13}.
\item[{$\Leftrightarrow$}] $A$ is an $(\gA^{\sur})^{\lup}$--bounded set
in $X$ by Theorems \ref{T:TB} and \ref{thm:4.15}.
\item[{$\Leftrightarrow$}]
$J_A\in[(\gA^{\sur})^{\lup}]^{\bsur}(X(A),X)$ by Theorem \ref{T:TB}.
\item[{$\Leftrightarrow$}] $J_A\in(\gA^{\lup})^{\bsur}(X(A),X)$ by
\cite[Proposition  3.5]{WW93}.
\item[{$\Leftrightarrow$}] $A$ is an $\gA^{\lup}$--bounded set in $X$
by Theorem \ref{T:TB}.
\end{enumerate}
\end{proof}





\begin{proposition}\label{prop:6.5} Let $\cM$ be a generating bornology on
LCS's and
$\gA=\sO^b(\cM)$. The $\gA$--topology coincides with the Grothendieck
topology
generated by $\gA$ on every LCS.
\end{proposition}
\begin{proof} It is easy to see that $\gA\subseteq\gA^{\rup}_{\bb}$. The
result
follows from Theorems \ref{thm:4.9} and \ref{thm:4.16}.
\end{proof}

\begin{proposition}\label{prop:6.6} Let $\cM$ be a generating bornology on
LCS's and
$\gA=\sO(\cM)$. Then the $\gA$--topology coincides with the
$(\gA^{\inj}_{\bb})^{\sup}$--topology on each infracomplete LCS.
\end{proposition}
\begin{proof} Since $\gA\subseteq(\gA^{\inj}_{\bb})^{\sup}$, $\sT(\gA)$ is
always
weaker than $\sT((\gA^{\inj}_{\bb})^{\sup})$ on each LCS. By
\cite[Corollary  3.2]{WW93},
$(\gA^{\inj}_{\bb})^{\sup}$ is injective. Let $p$ be a $\sT
((\gA^{\inj}_{\bb})^{\sup})$--continuous seminorm on an infracomplete
LCS $X$.
Then $\wQ_p\in(\gA^{\inj}_{\bb})^{\sup}(X,\wX_p)$. Let $B$ be a bounded
$\s$--disk in $X$. Now $J_{\wX_p}\wQ_pJ_B\in\gA^{\inj}_{\bb}(X(B),
\wX^{\inj}_p)$ implies $\wQ_pJ_B\in\gA^{\inj}_{\bb}(X(B),\wX_p)$ and
then again
implies $J_{\wX_p}\wQ_pJ_B\in\gA_{\bb}(X(B),
\wX^{\inj}_p)$ by Proposition \ref{hulls}.  Consequently, $J_{\wX_p}
\wQ_p\in\gA_{\bb}=\sO(\cM_{\bb})$. It turns out that
$\wQ_p\in
\gA^{\inj}$, or equivalently, $p$ is an $\gA$--continuous seminorm by
Theorem \ref{T:TB}.
\end{proof}

\begin{example}\label{eg:6.7.} Let $X=\bk^{(I)}$ be the locally convex direct sum
of card
$(I)$ many $\bk$'s where the index set $I$ is uncountable. $X$ is
infracomplete.
Let $\cM_{pc}$ be the generating bornology of precompact sets ($=$
totally bounded sets). Then $\sO^b
(\cM_{pc})=\gK_p$, the ideal of all precompact operators and
$\sO(\cM_{pc})=\gK^{\loc}_p$,
the ideal of all locally precompact operators, \ie those sending
bounded sets
onto precompact sets. $\gK_p$ is surjective but not
bornologically surjective and $\gK^{\loc}_p$ is bornologically
surjective. Now
$id_X\in\gK^{\loc}_p$ implies $X$ is a $\gK^{\loc}_p$--topological
space. On
the other hand, $X$ is not a $\gK_p$--topological space (cf.\
\cite[p.~40]{HN81}). This serves as a counter--example of
$\gA^{\sup}$--topology
$=\gA^{\inf}$--topology and $\gA^{\sup}$--topological spaces
$=\gA^{\inf}$--topological spaces,
although we always have
$\gA^{\rup}$--topology
$=\gA^{\inf}$--topology and $\gA^{\rup}$--topological spaces
$=\gA^{\inf}$--topological spaces.
By the way, $X$ is both $\gK_p$--bornological and $\gK^{\loc}_p$--
bornological, \ie a co--Schwartz space but not a Schwartz
space.
\end{example}

\begin{proposition}\label{prop:6.8} Let $\sG$ be a generating topology on
LCS's and
$\gA=\sO^b(\sG)$. The $\gA$--bornology coincides with the Grothendieck
bornology generated by $\gA$ on each infracomplete LCS.
\end{proposition}
\begin{proof} It follows from the easy fact $\gA\subseteq
\gA^{\lup}_{\bb}$ and Theorems \ref{thm:4.13} and \ref{thm:4.16}.
\end{proof}

\begin{proposition}\label{prop:6.9} Let $\sG$ be a generating topology on
LCS's and
$\gA=\sO(\sG)$. Then the $\gA$--bornology coincides with the
$(\gA^{\sur}_{\bb})^{\sup}$--bornology on
every infracomplete LCS.
\end{proposition}
\begin{proof} Similar to Proposition \ref{prop:6.6}.
\end{proof}

\begin{example}\label{eg:6.10.} Let $X=\bk^I$ be the product space of card$(I)$
many $\bk$'s
where the index set is uncountable. $X$ is infracomplete. Let
$\sG_{pc}$ be the
generating topology defined by the precompact seminorms,
\ie $\sG_{pc}=\sT(\gK_p)$, where $\gK_p$ is the ideal of all
precompact
operators between LCS's (see Wong
\cite{Wong79}). Then $\sO^b(\sG_{pc})$ is the ideal
$\gK^b_p$ of all
{\it quasi--Schwartz\/} ($=$ {\it precompact--bounded\/}, cf.\ Rankte
\cite{Ran72}) {\it operators\/} between LCS's. $\sO(\sG_{pc})$ is the
ideal of
those continuous operators between LCS's which are still continuous
when the
domain space $X$ equipped with the (coarser) precompact topology
$\sG_{pc}(X)$.
$X$ is not a $\gK^b_p$--bornological space since otherwise (by Theorem
\ref{thm:3.10}) we
would have the canonical embedding from $\bk^{(I)}$ into $\bk^I$ being
quasi--Schwartz and this is not the case as shown in \cite[p.\ 399]{Jun83}. $X$
is, however, an $\sO(\sG_{pc})$--bornological space since all bounded
sets in
$X$ are precompact. This serves as a counter--example of
$\gA^{\sup}$--bornology $=\gA^{\inf}$--bornology and $\gA^{\sup}$--
bornological
spaces $=\gA^{\inf}$--bornological spaces,
although we always have $\gA^{\rup}$--bornology $=\gA^{\inf}$--
bornology on every infracomplete LCS.
By the way, $X$ is both
$\gK^b_p$--topological and $\sO(\sG_{pc})$--topological, \ie a
Schwartz space
but not a co--Schwartz space.
\end{example}

\begin{remark} It may be interesting to study the $\gA^{\sup}$--topology
and the
$\gA^{\sup}$--bornology for an operator ideal $\gA$ on Banach spaces.
Propositions \ref{prop:6.6} and \ref{prop:6.9} suggest the conjectures that $\sO(\sT(\gA))=
(\gA^{\inj}_{\bb})^{\sup}$ and
$\sO(\sB(\gA))=(\gA^{\sur}_{\bb})^{\sup}$ where
$\gA$ is an operator ideal on LCS's.
\end{remark}

\begin{theorem}\label{thm:5.1} Let $\gA$ be an operator ideal on LCS's, and
$X$ a
LCS. Then
\begin{enumerate}[(a)]
\item  $X$ is $\gA^{\dual}$--bornological $\Rightarrow X'_\b$ is
$\gA$--topological.
\end{enumerate}
\noindent If, in addition, $X$ is infrabarrelled then
\begin{enumerate}[(b)]
\item  $X$ is $\gA^{\dual}$--topological $\Rightarrow X'_\b$ is
$\gA$--bornological.
\end{enumerate}
\begin{enumerate}[(c)]
\item  $X'_\b$ is $\gA^{\dual}$--bornological $\Rightarrow X$ is
$\gA$--topological.
\end{enumerate}
\noindent If, in addition to all above, $\gA$ is also injective then
\begin{enumerate}[(d)]
\item  $X'_\b$ is $\gA^{\dual}$--topological $\Rightarrow X$ is
$\gA$--bornological
\end{enumerate}
\end{theorem}
\begin{proof}

(a) Let $V$ be an absolutely convex, closed $0$--neighborhood in
$X'_\b$. Then $V^\circ$ is bounded and hence $\gA^{\dual}$--bounded in
$X$. So
that there is a normed space $N$, a $T$ in $\gA^{\dual}(N,Y)$ with
$TU_N\supseteq V^\circ$. Consequently, $(TU_N)^\circ=(T')^{-1}U_{N'}\subseteq
V^{\circ\circ}=V$ and $T'\in\gA(X'_\b,N')$. It follows that $V$ is an
$\gA$--neighborhood of $0$ in $X'_\b$.

(b) Let $B$ be a bounded set in $X'_\b$. Then $B^\circ$ is a closed
bornivorous
barrel in $X$, and hence a $0$--neighborhood, and consequently an
$\gA^{\dual}$--neighborhood of $0$ in $X$. Therefore there is a Banach
space
$F$, an $T$ in $\gA^{\dual}(X,F)$ such that $B^\circ\supseteq T^{-1}(U_{F})$. It
follows $B\subseteq B^{\circ\circ}\subseteq T'U_{F'}$. Since $T'\in\gA
(F',X'_\b)$, $B$ is $\gA$--bounded in $X'_\beta$.

(c) Let $V$ be an absolutely convex, closed $0$--neighborhood in $X$.
Then $V^\circ$ is bounded and hence $\gA^{\dual}$--bounded in $X'_\b$.
So that
there is a normed space $N$ and a $T$ in $\gA^{\dual}(N,X'_\b)$ such
that
$TU_N\supseteq V^\circ$ and thus $(T')^{-1}U_{N'}\subseteq
V^{\circ\bullet}$,
where $V^{\circ\bullet}$ is the polar of $V^\circ$ in the strong bidual
$X''_{\b\b}$ of $X$. It follows that $V^{\circ\bullet}$ is an
$\gA$--neighborhood of $0$ in $X''_{\b\b}$. Since the evaluation map
$K_X
:X\to X''_{\b\b}$ is continuous, $K^{-1}_X(V^{\circ\bullet})=V$ is an
$\gA$--neighborhood of $0$ in $X$ by (GT$_2$).

(d) Let $B$ be an absolutely convex bounded set in $X$. It suffices to
check that $J_B\in\gA^{\bsur}(X(B),X)$. Note that $B^\circ$ is a
neighborhood of $0$, and hence an $\gA^{\dual}$--neighbor\-hood of $0$ in
$X'_\b$.
Hence there exist a Banach space $F$ and a $T$ in $\gA^{\dual}(X'_\b,F)$
such that
$B^\circ\supseteq T^{-1}(U_{F})$. Hence $B^{\circ\bullet}\subseteq
T'U_{F'}$.
Since $T'\in\gA(F',X''_{\b\b})$,
$B^{\circ\bullet}$ is $\gA$--bounded in
$X''_{\b\b}$ and thus  $J_{B^{\circ\bullet}}\in\gA^{\bsur}(X''_{\b\b}
(B^{\circ\bullet}),X''_{\b\b})$. Now $K_XB\subseteq B^{\circ\bullet}$
ensures
that there is a $K_B$ in $\sL(X(B),X''_{\b\b} (B^{\circ\bullet}))$
such that $J_{B^{\circ\bullet}}K_B=K_XJ_B$. Hence
$K_XJ_B\in\gA^{\bsur}(X(B),
X''_{\b\b})$ and it follows
$$
J_B\in(\gA^{\bsur})^{\inj}(X(B),X)
=(\gA^{\inj})^{\bsur}(X(B),X)=\gA^{\bsur}(X(B),X),
$$
by Proposition \ref{hulls}.
\end{proof}

\begin{theorem}\label{thm:5.2} Let $\gA$ be a symmetric operator ideal (\ie
$\gA
\subseteq\gA^{\dual}$) on LCS's, and $X$ an infrabarrelled LCS. Then
\begin{enumerate}[(a)]
\item  $X$ is $\gA$--topological $\Leftrightarrow X'_\b$ is
$\gA$--bornological.
\end{enumerate}
\noindent If, in addition, $\gA$ is injective, then
\begin{enumerate}[(b)]
\item  $X$ is $\gA$--bornological $\Leftrightarrow X'_\b$ is
$\gA$--topological.
\end{enumerate}
\end{theorem}
\begin{proof} A consequence of Theorem \ref{thm:5.1}.
\end{proof}

\begin{theorem}\label{thm:6.14} Let $\sG$ be a generating topology on LCS's
and $X$
be an infrabarrelled LCS. Then
\begin{enumerate}[(a)]
\item $X$ is $\sO^b(\sG^\circ)$--topological $\Rightarrow
X'_\b$ is
$\sO^b(\sG)$--bornological.
\item $X$ is $\sO^b(\sG^\circ)$--bornological $\Rightarrow
X'_\b$
is $\sO^b(\sG)$--topological.
\item $X'_\b$ is $\sO^b(\sG^\circ)$--topological
$\Rightarrow X$ is
$\sO^b(\sG)$--bornological.
\item $X'_\b$ is $\sO^b(\sG^\circ)$--bornological
$\Rightarrow X$
is $\sO^b(\sG)$--topological.
\end{enumerate}
\noindent
In case $\sO(\sG)$ is symmetric or $\sG$ has the subspace property, all
above
implications become equivalences.
\end{theorem}
\begin{proof} We prove (c) only and all others are similar. Suppose
$X'_\b$ is $\sO^b(\sG^\circ)$--topological. By Proposition \ref{P:TBdual}(a), $X'_\b$ is also
$\sO^b
(\sG)^{\dual}$--topological. Note that $\sO^b(\sG)$ is injective. Hence
by
Theorem \ref{thm:5.1}(d), $X$ is $\sO^b(\sG)$--bornological. In case $\sO(\sG)$
is
symmetric or $\sG$ has the subspace property, if $X$ is $\sO^b
(\sG)$--bornological, $X$ is also $\sO^b(\sG^\circ)^{\dual}$--
bornological by Proposition \ref{P:TBdual}(b)$'$. By
Theorem \ref{thm:5.1}(a), $X'_\b$ is $\sO^b(\sG^\circ)$--topological, as
asserted.
\end{proof}

\begin{example}\label{eg:6.15.} Let $\sG_\s$ be the generating system of
$\s(X,X')$--topology on each LCS $X$. $\sG^\circ_\s(X)$ is thus the
convex
bornology $\cM(X)$ consisting of those bounded subsets $B$ of $X$ whose
polars
$B^\circ$ are $\s(X'_\b,X''_{\b\b})$--neighborhoods of $0$ in $X'_\b$.
Now
both $\sO^b(\sG_\s)$ and $\sO^b(\sG^\circ_\s)$ define the ideal $\gF$
of continuous
operators of finite rank. Moreover, $\sG_\s$ has the subspace property.
Theorem \ref{thm:6.14} applies and says that for an infrabarrelled LCS $X$, we have
\begin{enumerate}[(a)]
\item $X$ is $\gF$--topological if and only if $X'_\b$ is
$\gF$--bornological;
\item $X$ is $\gF$--bornological if and only if $X'_\b$ is
$\gF$--topological.
\end{enumerate}
Unlike the case of Banach spaces, an $\gF$--topological or
$\gF$--bornological
LCS need not be of finite dimension. For examples, the LCS $\bk^I$
is
$\gF$--topological and the LCS $\bk^{(I)}$ is
$\gF$--bornological. This is because the weak topology
$\s(\bk^I,\bk^{(I)})$
of $\bk^I$ coincides with the product topology of $\bk^I$ and every
bounded set
in $\bk^{(I)}$ is of finite dimension. Here the index set $I$ is
arbitrary.
\end{example}


\section{Examples and Applications}\label{s:egs}

This last section is devoted to  examples and
applications, showing the powerful techniques developed in the previous
sections.  Many other elegant applications of the theory of
Grothendieck spaces and co--Grothendieck
spaces can be found, for example, in \cite{Jar81}, \cite{Jun83}, \cite{Lin88}, \cite{Pie72},
\cite{Pie80}. The concepts  of $\gA$--topological spaces
and
$\gA$--bornological spaces are also well--developed in the context in
\cite{HN81},
\cite{Pie72}, \cite{Ran72}, \cite{Wong76}, \cite{WW88}, and, in
particular,
\cite{Wong79}.


\subsection{Schwartz spaces and co-Schwartz spaces}

\begin{definition}[{see, e.g., \cite[{p.~14}]{Wong79}}] A continuous seminorm $p$
on a
LCS $X$ is said to be {\it precompact\/} if there exists a $(\lam_n)$
in $c_0$
and an equicontinuous sequence $\{x'_n\}$ in $X'$ such that
$$
p(x)\leq\sup\{|\lam_n\langle x,x'_n\rangle| :n\geq1\},\qquad\forall
x\in X.
$$
\end{definition}

Denote by $\sG_{pc}(X)$ the locally convex (Hausdorff) topology on $X$
defined
by all precompact seminorms on $X$. It is easy to see that
$\sG_{pc}=\{\sG_{pc}
(X) :X$ is a $LCS\}$ is a generating topology. It is a classical result
(cf.\ \cite{Ran72} or \cite{Wong79}) that $p$ is a precompact seminorm on a
LCS $X$
if and only if the canonical map $Q_p :X\to X_p$ is precompact.
Then, $\gK_p=\sO(\sG_{pc})$ is the ideal of all precompact operators, and
$\gK^b_p=\sO^b(\sT(\gK_p))$ is
the ideal of all quasi--Schwartz (\ie precompact--bounded) operators
between
LCS's.

\begin{definition} A LCS $X$ is said to be a {\it Schwartz space\/} if
every
continuous seminorm $p$ on $X$ is precompact.
\end{definition}

     We provide a new proof of the following classical result.

\begin{theorem}[{see \cite[ pp.\ 17 and 26]{Wong79}}]\label{thm:7.1}
Let $X$ be a
LCS. The
following are all equivalent.
\begin{enumerate}[(a)]
\item  $X$ is a Schwartz space.
\item  For each continuous seminorm $p$ on $X$ there is a
continuous
seminorm $q$ on $X$ such that $p\leq q$ and the canonical map $Q_{pq}$
belongs
to $\gK_p(X_q,X_p)$.
\item  $Q_p\in\gK_p(X,X_p)$ for every continuous seminorm
$p$ on $X$.
\item  For any $0$--neighborhood $U$ in $X$ there exists a
$0$--neighborhood $V$ in $X$ such that $V\subseteq U$ and the canonical
map from
$X'(U^\circ)$ into $X'(V^\circ)$ is precompact.
\item  $\sL(X,N)=\gK_p(X,N)$ for every normed (or Banach)
space $N$.
\item  $\gK^b_p(X,Y)=\sL^b(X,Y)$ for every LCS $Y$.
\item  $\gK^b_p(X,Y)=\gK_p(X,Y)$ for every LCS $Y$.
\item  $\sL(X,N)=\gK^b_p(X,N)$ for every normed (or Banach)
space
$N$.
\item  $X$ is a $\gK_p$--topological space.
\end{enumerate}
\end{theorem}
\begin{proof}
(a)$\Leftrightarrow$(c)$\Leftrightarrow$(e)$\Leftrightarrow$(i) are
due to Theorem \ref{thm:3.4} and the injectivity of $\gK_p$.
(a)$\Leftrightarrow$(b)
because $\gK_p=\sO^b(\cM_{pc})$ where $\cM_{pc}$ is the generating
bornology of
precompact sets and Proposition \ref{prop:6.5} applies. (b)$\Leftrightarrow$(d)
follows
from the complete symmetry of the restriction $(\gK_p)_{\bb}$ of
$\gK_p$ to
Banach spaces, \ie $(\gK_p)^{\dual}_{\bb}=(\gK_p)_{\bb}$.
(a)$\Leftrightarrow$(f)$\Leftrightarrow$(h) are consequences of Theorem
\ref{thm:3.4}.
(i)$\Rightarrow$(g) is contained in Proposition \ref{P:OTOOBO} and Theorem  \ref{thm:6.4}. Finally, for
(g)$\Rightarrow$(h),
denote by
$N_\s$ the LCS $(N,\s(N,N'))$. For every $T$ in $\sL(X,N)$,
$T\in\sL(X,N_\s)=\gK_p
(X,N_\s)$. Hence, by (g), $T\in\gK^b_p(X,N_\s)=\gK^b_p(X,N)$ since $N$
and
$N_\s$ carry the same (von Neumann) bornology.
\end{proof}

\begin{definition} A LCS $Y$ is said to be a {\it co--Schwartz space\/}
if its
strong dual $Y'_\b$ is a Schwartz space.
\end{definition}

\begin{theorem}\label{thm:7.2} Let $Y$ be a LCS. Consider the following
statements.
\begin{enumerate}[(a)]
\item  $Y$ is a co--Schwartz space.
\item For each bounded disk $B$ in $Y$ there is a bounded
disk $A$
in $Y$ with $B\subseteq A$ such that the canonical map $J_{AB}$ from
$Y(B)$ into
$Y(A)$ belongs to $\gK_p(Y(B),Y(A))$.
\item  $J_B\in\gK_p(Y(B),Y)$ for each bounded disk $B$ in
$Y$.
\item  $\sL(N,Y)=\gK_p(N,Y)$ for every normed space $N$.
\item  $\sL^b(X,Y)=\gK_p(X,Y)$ for every LCS $X$.
\item  $Y$ is a $\gK_p$--bornological space.
\end{enumerate}
\noindent
We have (a)$\Leftrightarrow$(b)$\Rightarrow$(c)$\Leftrightarrow$
(d)$\Leftrightarrow$(e)$\Leftrightarrow$(f).
\end{theorem}
\begin{proof} (a)$\Leftrightarrow$(b) follows from the equivalence
(a)$\Leftrightarrow$(b) in the last theorem and the complete symmetry
of $(\gK_p)_{\bb}$.
(c)$\Leftrightarrow$(d)$\Leftrightarrow$(e)$\Leftrightarrow$(f)
are just  examples of Theorem \ref{thm:3.10}. (b)$\Rightarrow$(c) is trivial.
Finally, the
LCS $\bk^I$, where the index set $I$ is uncountable, furnishes a
counter--example of the missing implication.
\end{proof}

\begin{proposition}\label{prop:7.3} Let $X$ be an infrabarrelled LCS. Then $X$
is a
Schwartz space ($\resp$ co--Schwartz space) if and only if $X'_\b$ is a
co--Schwartz space ($\resp$ Schwartz space) if and only if $X''_{\b\b}$
is a
Schwartz space ($\resp$ co--Schwartz space).
\end{proposition}
\begin{proof} Repeat applying Theorems \ref{thm:7.1} and \ref{thm:7.2}, and the complete
symmetry of $(\gK_p)_{\bb}$.
\end{proof}

\begin{remark} Besides the ideals of precompact operators and
quasi--Schwartz
operators, one can also employ the ideal $\sL_{\text{im}}$ of {\it limit
operators\/} to
define Schwartz spaces and co--Schwartz spaces. See \cite{Lin88} for some other
internal characterization of Schwartz spaces due to the introduction of
$\sL_{\text{im}}$.
\end{remark}

Similar to Schwartz space we can relate infra--Schwartz spaces to the
ideal
$\gW$ of weakly compact operators between LCS's. Incidentally, readers
should
have no difficulty to figure out that $\gK_p$--bornological spaces are, in fact,
semi--Montel spaces and $\gW$--bornological spaces are exactly
semi--reflexive
spaces. We leave these to the interested readers and refer them to
\cite{HN81} for
more information about the classical theory of these spaces.

\subsection{Nuclear spaces and co-nuclear spaces}

\begin{definition} A continuous seminorm $p$ on a LCS $X$ is called an
{\it absolutely summing seminorm\/} ($=$ {\it prenuclear seminorm\/} in
\cite{Wong79}) if there exists a $\s(X',X)$--closed equicontinuous subset
$B$ of $X'$ and a positive Radon measure $\mu$ on $B$ such that
$$
p(x)\leq\int_B|\langle x,x'\rangle|d\mu(x'),\qquad\forall x\in X.
$$

Let $\sG_{\as}(X)$   be the locally convex (Hausdorff)
topology on $X$ generated by the family of all absolutely summing
seminorms on
$X$. It is easy to see that the system $\sG_{\as}=\{\sG_{\as}(X) :X$ a
$LCS\}$
is a generating topology. A continuous operator $T$ from a LCS $X$
into a
LCS $Y$ is said to be {\it absolutely summing\/} if
$T\in\sO(\sG_{\as})
(X,Y)=\sL(X_{\sG_{\as}},Y)$. In case $X$ and $Y$ are Banach spaces, $T$
is
absolutely summing if and only if $T$ sends every weakly summable
series in $X$
to an absolutely summable series in $Y$. Denote by $\gP=
\sO(\sG_{\as})$ the
injective ideal
of all absolutely summing operators between LCS's, and
by $\gP^b=\sO^b(\sG_{\as})$ the injective ideal  of
{\it prenuclear--bounded operators\/} \cite{Wong79}.

A continuous operator $T$ from a LCS $X$ into a LCS $Y$ is said to be {\it nuclear\/} if there exist a $(\lam_n)$ in
$\l_1$, an
equicontinuous sequence $\{a_n\}$ in $X'$ and a sequence $\{y_n\}$
contained in
an infracomplete bounded disk $B$ in $Y$ such that
$T=\Sigma_{n}\lam_na_n\otimes y_n$, \ie
$Tx=\Sigma_{n}\lam_na_n(x)y_n$ for each $x$ in $X$.
Denote by
$\gN$ the ideal of all nuclear operators between LCS's. Note that
$\gN$ is
symmetric. It is more or less classical that $\gP=
\gP^{\rup}_{\bb}$ \cite[ p.\ 76]{Wong79}, $\gN=\gN^{\inf}_{\bb}$
\cite[p.\ 144]{Wong79} and $\gP^3_{\bb}\subset \gN_{\bb}\subset \gP_{\bb}$
\cite[p.\ 145]{Wong79} (in fact, we have $\gP^2_{\bb}\subset \gN_{\bb},
(\gN^{\inj}_{\bb})^2\subset \gN_{\bb}$, cf.\ \cite{Pie72}).
\end{definition}

\begin{definition} A LCS $X$ is said to be {\it nuclear\/} if every
continuous
seminorm $p$ on $X$ is absolutely summing. A LCS $Y$ is said to be
{\it co--nuclear\/} if its strong dual $Y'_\b$ is nuclear.
\end{definition}

     We provide a new proof of the following classical result.

\begin{theorem}[{see \rm{\cite[pp.\ 149 and 157]{Wong79}}}]\label{thm:7.4}
Let $X$ be a
LCS. The
following are all equivalent.
\begin{enumerate}[(a)]
\item  $X$ is a nuclear space.
\item  $Q_p\in\gP(X,X_p)$ for every continuous seminorm $p$
on $X$.
\item  For each continuous seminorm $p$ on $X$ there exists
a
continuous seminorm $q$ on $X$ with $p\leq q$ such that the canonical
map
$Q_{pq}\in\gP(X_q,X_p)$.
\item  $id_X\in\gP(X,X)$.
\item  $\gP(X,Y)=\sL(X,Y)$ for every LCS $Y$.
\item  $\gP(X,N)=\sL(X,N)$ for every normed space $N$.
\item  $\wQ_p\in\gN(X,\wX_p)$ for every continuous seminorm
$p$ on $X$.
\item  $\gN(X,F)=\sL(X,F)$ for every Banach space $F$.
\item  For each continuous seminorm $p$ on $X$ there exists
a
continuous seminorm $q$ on $X$ with $p\leq q$ such that the canonical
map
$\wQ_{pq}\in\gN(\wX_q,\wX_p)$.
\item  For each $0$--neighborhood $V$ in $X$ there is a
$0$--neighborhood $U$ in $X$ with $U\subseteq V$ such that the canonical
map
$X'(V^{\circ})\to X'(U^{\circ})$ is nuclear.
\item  $\gP^b(X,Y)=\sL^b(X,Y)$ for every LCS $Y$.
\item $\sL^b(X,Y)\subseteq\gP(X,Y)$ for every LCS $Y$.
\item  $\gK_p(X,Y)\subseteq\gP^b(X,Y)$ for every LCS $Y$.
\item  $\sL^b(X,Y)\cap\gP(X,Y)=\gP^b(X,Y)$ for every LCS $Y$.
\item  $X$ is a $\gP$--topological space.
\item  $X$ is a $\gN$--topological space.
\end{enumerate}
\end{theorem}
\begin{proof} (a) $\Leftrightarrow$ (b) $\Leftrightarrow$ (d)
$\Leftrightarrow$
(e) $\Leftrightarrow$ (f) $\Leftrightarrow$ (l) $\Leftrightarrow$ (o)
$\Leftrightarrow$ (p)
are due to Theorem \ref{thm:3.4}. Since
$\gP=\gP^{\rup}_{\bb}$, we have (a) $\Leftrightarrow$ (c) by Theorem
\ref{thm:4.16}.  (a) $\Leftrightarrow$ (n) $\Leftrightarrow$ (k) are due to
Proposition \ref{P:OTOOBO} and Theorem \ref{thm:6.4}. (k) $\Rightarrow$ (m) is obvious. To prove (m)
$\Rightarrow$ (k)
we employ the same trick as in Theorem \ref{thm:7.1}. (c) $\Leftrightarrow$ (i)
follows from the fact that
$\gP^3\subset \gN\subset \gP$. (i) $\Rightarrow$ (g) $\Rightarrow$ (h) are trivial.
(h)
$\Rightarrow$ (l) because $\gN\subseteq\gP$ and $\gP$ is injective. (i)
$\Rightarrow$ (j) is ensured by the symmetry of $\gN_{\bb}$.

Finally,
we prove
(j) $\Rightarrow$ (i). Let $V_p=\{x\in X :p(x)\leq1\}$ be the
$0$--neighborhood
associated to a continuous seminorm $p$ on $X$. By (j), there is a
continuous
seminorm $q$ on $X$ such that $V_q\subseteq V_p$ (\ie $p\leq q$) and
$\wQ'_{pq}
:X'(V^\circ_p)\to X'(V^\circ_q)$ is nuclear. By the symmetry of
$\gN_{\bb}$,
$\wQ''_{pq} :(X'(V^\circ_q))'\to(X'(V^\circ_p))'$ is nuclear, too.
Hence
$\wQ''_{pq}$ is absolutely summing. Now $(X'(V^\circ_q))'$ and
$(X'(V^\circ_p))'$ are isometrically isomorphic to $X''_q$ and $X''_p$,
respectively. By the injectivity of $\gP$, $\wQ_{pq}$ is absolutely
summing.
Repeating the same argument, we shall have continuous seminorms $q_1$
and $q_2$
on $X$ such that $q\leq q_1\leq q_2$ and $\wQ_{qq_1}$ and
$\wQ_{q_1q_2}$ are
both absolutely summing. Now $p\leq q_2$ and
$\wQ_{pq_2}=\wQ_{pq}\wQ_{qq_1}
\wQ_{q_1q_2}\in\gP^3_{\bb}\subseteq\gN_{\bb}$, and we are
done.
\end{proof}

\begin{remark} There are concepts of {\it quasi--nuclear--seminorms,
quasi--nuclear operators\/} and {\it quasi--nuclear--bounded
operators,\/}
cf.\ \cite{Wong79}. They can be used to define nuclear spaces like $\gP$ and $\gN$.
However,
they are simply, respectively, the $\gN$--seminorms, $\gN^{\inj}$--operators and
$(\sT(\gN^{\inj}))^b$--operators. Using the same kind of argument in
Theorem \ref{thm:7.4},
one can easily prepare a longer list of equivalences. We leave this
 to the interested readers.
\end{remark}

\begin{theorem}\label{thm:7.5} Let $Y$ be an infrabarrelled LCS. The
following are
all equivalent.
\begin{enumerate}[(a)]
\item  $Y$ is a co--nuclear space.
\item  For each bounded disk $B$ in $Y$ there is a bounded
disk $A$
in $Y$ with $B\subseteq A$ such that the canonical map $J_{AB}$ from
$Y(B)$ into
$Y(A)$ is nuclear.
\item  $J_B\in\gN(Y(B),Y)$ for every bounded disk $B$ in
$Y$.
\item  $\sL(N,Y)=\gN(N,Y)$ for every normed space $N$.
\item  $\sL^b(X,Y)\subseteq\gN(X,Y)$ for every LCS $X$.
\item  $Y$ is an $\gN$--bornological space.
\end{enumerate}
\end{theorem}
\begin{proof} Assume first that $Y$ is co--nuclear and $B$ is a bounded
disk in
$Y$. Then $B^\circ$ is a $0$--neighborhood in $Y'_\b$. Hence there is a
bounded
disk $A$ in $Y$ with $B\subseteq A$ such that the canonical map $Y''
(B^{\circ\bullet}) \to Y''(A^{\circ\bullet})$ is absolutely summing.
Since $\gP$ is
injective, the canonical map $J_{AB}$ from $Y(B)$ into $Y(A)$ is also
absolutely summing. Do this twice more and we shall get (a)
$\Rightarrow$ (b)
since $\gP^3_{\bb}\subseteq\gN_{\bb}$. (b) $\Rightarrow$ (c), (d), (e)
and each
one of them $\Rightarrow$ (f) are straightforward. We consider (f)
$\Rightarrow$ (a). Note that $\gN$ is symmetric. Now Theorem \ref{thm:5.1}(a)
gives the
desired conclusion.
\end{proof}

\begin{proposition}\label{prop:7.6} Let $X$ be an infrabarrelled LCS. Then $X$
is a
nuclear space ($\resp$ co--nuclear space) if and only if $X'_\b$ is a
co--nuclear space ($\resp$ nuclear space) if and only if $X''_{\b\b}$
is a
nuclear space ($\resp$ co--nuclear space).
\end{proposition}
\begin{proof} In view of Theorem \ref{thm:6.14}, it suffices to mention that the
generating system $\sG_{\as}$ of absolutely summing topology has the
subspace
property. As a result, an infrabarrelled LCS $X$ is nuclear if and
only if
$X''_{\b\b}$ is nuclear. The other implications follow from
this.
\end{proof}

Since $\gN\subset \gK_p$ we have the well--known

\begin{proposition}\label{prop:7.7} All nuclear ($\resp$ co--nuclear) spaces
are Schwartz ($\resp$ co--Schwartz) spaces.
\end{proposition}

\subsection{Permanence properties}

We collect some results from \cite{Jun83} about the permanence properties of
Grothendieck spaces and co--Grothendieck spaces.

\begin{theorem}[{\cite[Junek]{Jun83}}]\label{thm:7.8}
Let $\gA$ be an operator ideal
on
Banach spaces.
\begin{enumerate}[(a)]
\item Any product of $\Groth(\gA)$--spaces is a $\Groth(\gA)$--space.
\item Any locally convex direct sum of
co--$\Groth(\gA)$--spaces is
a co--$\Groth(\gA)$--space.
\item If $\gA$ is equivalent to some injective ideal then
any
subspace of a $\Groth(\gA)$--space is a $\Groth(\gA)$--space.
\item If $\gA$ is equivalent to some surjective ideal then
any
quotient space of a co--$\Groth(\gA)$--space is a
co--$\Groth(\gA)$--space.
\item If $\gA$ is injective then any projective limit of
$\Groth(\gA)$--spaces is a $\Groth(\gA)$--space.
\item If $\gA$ is injective then any subspace of a
co--$\Groth(\gA)$--space is a co--$\Groth(\gA)$--space.
\end{enumerate}
\end{theorem}

If one applies them together with other results in the earlier parts of
this paper to Schwartz spaces, infra--Schwartz spaces, nuclear spaces,
and their ``co--spaces", one can obtain a long list of permanence
properties
of these spaces, cf. \cite{Jar81} or \cite{Jun83}.

\subsection{Other applications}

Along the same line of reasoning in this paper one can develop similar
applications of operator ideals to the theory of tensor products,
partially ordered locally convex spaces and $C^*$--algebras.

It is of no doubt that the initial idea of operator ideals comes from
tensor
products. In \cite[p.\ 49]{Mic78}, Michor suggested a method to
construct a
tensor norm $\gA^{\otimes}$ associated to each operator ideal $\gA$ on
Banach
spaces (see \cite{Pie80} for details about
{\it quasi--normed ideal\/}).  See also
\cite{Wong79} and  those famous works of
A.\ Grothendieck and R. Schatten.

Let $E$ and $F$ be Banach
spaces and $\gA$ be an operator ideal on Banach spaces with ideal norm
$\a$.
Define $\V\cdot\V_{\gA^\otimes}$ on $E\otimes F$ by
$$
\V\sum x_i\otimes y_i\V_{\gA^\otimes}=\sup\{|\sum\langle
y_i,Tx_i\rangle|
:T\in\gA(E,F'),\quad\a(T)\leq1\}.
$$
$\V\cdot\V_{\gA^\otimes}$ turns out to be a reasonable cross norm. We
denote by $E{\otimes}_{\gA}F$ the $\gA$--{\it tensor product\/} of $E$
and $F$, that is, the completion of $E\otimes F$ under
$\V\cdot\V_{\gA^\otimes}$.  Y.~C.~Wong
\cite[p.\ 279]{Wong79} showed that if $\gA$ is the normed ideal $(\gP,P)$
of all
absolutely summing operators between Banach spaces, we would have
$(E\otimes_{\gP}F)'\cong(\gP(E,F'),P)$.

In general, let $\gA$ be an
operator
ideal on LCS's. We can define a tensor product topology associated to
$\gA$
by a family of $\gA$--bilinear forms. A continuous bilinear form $b$
 on $X\times Y$ is said to be  an
$\gA$--{\it bilinear form,\/} if there is a $T$ in $\gA(X,Y'_\b)$ such
that
$b(x,y)=\langle y,Tx\rangle$. We write $b=b_T$ in this case. Detailed
properties of $b_T$ can be found in \cite{Wong79}. See also \cite{Jam87} for
other
comments. If $\gA$ is equipped with some locally convex topology (see
\cite{Jun83})
then we can define similar seminorms like the one as
$\V\cdot\V_{\gA^\otimes}$.
It might be interesting to investigate this kind of theory.

There is also an established theory of ideal topologies on partially
ordered
locally convex spaces. We give only one
example here and refer interested readers to \cite{Wong76}. Let
$(X,X_+,\sT)$
be a locally solid space.  A continuous seminorm $p$ on $X$ is said to
be a
$(PL)$--{\it seminorm\/} if there exists a positive $f$ in $X'$ such
that
$p(x)\leq\sup\{g(x) :-f\leq g\leq f\}$, $\forall x\in X$. It turns out
that a
continuous seminorm $p$ on a locally solid space $X$ is a
$(PL)$--seminorm if
and only if $Q_p$ is a cone--absolutely summing operator from $X$ onto
$X_p$.  Moreover, we have a list of characterizations of $\sT$ to be
the topology of uniform convergence on all order intervals as those
appeared in Theorems \ref{thm:7.1} and \ref{thm:7.4} (see \cite[p.~136]{Wong76}).
We would like to mention that in the case of partially ordered locally
convex
spaces, or Banach lattices, the correct concept of operator ideals may
be
the so--called {\it operator modules.\/} For more information about
operator
modules, see Schwarz \cite{Schwarz84}.

Finally, we finish this paper with a result of Jarchow \cite{Jar86}. Let $H$
be a Hilbert space and $A$ be a $C^*$--subalgebra of $B(H)$.

\begin{proposition}[{\cite[Jarchow]{Jar86}}]\label{prop:7.9}
The $\gW$--topology of
$A$,
\ie the ideal topology on $A$ generated by weakly compact operators,
is the
finest locally convex topology on $A$ which coincides with the
strong$^*$
(\ie the double strong) operator topology on bounded subsets of $A$.
The
completion of $A$ under this topology is
$(A^{**},\tau(A^{**},A^*))$.
\end{proposition}


\end{document}